\documentstyle[11pt, amssymb,amsmath,amstext,amsfonts,amscd,xypic]{article}

\makeatletter
\@addtoreset{equation}{section}
\makeatother

\def\e{{\mathsf e}}
\def\Gw{{\mathsf G}_W}
\def\bmu{{\boldsymbol\mu}_m}
\def\bchi{{\boldsymbol\chi}}
\def\po{\frac{1}{p}\cd }
\def\hp{\hphantom{x}}
\newcommand{\nc}{\newcommand}
\newcommand{\qed}{$\enspace\square$}

\newcommand{\dis}{\displaystyle}

\newcommand{\DR}{{\mathcal{DR}}}

\newtheorem{theorem}[equation]{Theorem}
\newtheorem{proposition}[equation]{Proposition}
\newtheorem{corollary}[equation]{Corollary}

\newtheorem{lemma}[equation]{Lemma}
\newtheorem{example}[equation]{Example}
\newtheorem{definition}[equation]{Definition}

\newenvironment{case}{\left\{\begin{array}{rl}}{\end{array}\right.}

\nc{\dlim}{\mathop{\lim\limits_{\longrightarrow}}}
\nc{\dlimk}{\mathop{\lim\limits_{{\longrightarrow \atop k}}}}
\nc{\dlimkl}{\mathop{\lim\limits_{{\longrightarrow \atop {k,l}}}}}
\nc{\FB}{{\mathcal F}{\mathsf D}}
\nc{\B}{{\mathsf D}}
\nc{\BS}{{\underline{\B}}}
\nc{\D}{{\mathcal D}}
\nc{\Brat}{{{\mathsf D}_{\mathrm frac}}}
\nc{\Grad}[1]{\mathsf{Gr}^{\mathsf{ad}}_{#1}}
\nc{\vdim}{{\mathrm v.}\dim}
\nc{\Qu}{{\mathsf Q}}
\nc{\FM}{{\mathfrak M}}
\nc{\TM}{{\widetilde M}}
\nc{\PPP}{{\PP^1\times_\tau\PP^1}}
\nc{\iz}{i_z}
\nc{\iw}{i_w}
\nc{\Quot}{{\mathrm Quot}^{\PPP}}
\nc{\BQuot}{{\mathrm Quot}^{\B}}
\nc{\Subs}{{{\mathrm Sub}_\CO^0}}
\nc{\BSubp}{{{\mathrm Sub}_\B^p}}
\nc{\BSubs}{{{\mathrm Sub}_\B}}
\nc{\CSubs}{{{\mathrm Sub}_\C}}
\nc{\CM}{{\mathcal M}}
\nc{\CU}{{\mathcal U}}
\nc{\epss}{{\eps^{-1}}}
\nc{\CG}{{\C\Gamma}}
\nc{\Diff}{{\mathop{\mathrm Diff}_\tau}}
\nc{\longleftrightarrows}{\quad\rlap{\raisebox{-2pt}{$\longleftarrow$}}%
\raisebox{2pt}{$\longrightarrow$}\quad}
\nc{\SS}{{\mathop{\sf S}}}

\def\beq{\begin{equation}}
\def\eeq{\end{equation}}

\newcommand{\iso}{{\;\stackrel{_\sim}{\longrightarrow}\;}}
\newcommand{\cd}{\!\cdot\!}
\newcommand{\vi}{${\sf {(i)}}\;$}
\newcommand{\vii}{${\sf {(ii)}}\;$}
\newcommand{\viii}{${\sf {(iii)}}\;$}
\newcommand{\id}{{\mbox{\sl Id}}}

\newcommand{\into}{\,\hookrightarrow\,}
\newcommand{\too}{\,\longrightarrow\,}
\newcommand{\onto}{\,\twoheadrightarrow\,}

\newcommand{\Hom}{{\mbox{\sl Hom}}}

\newcommand{\op}{\mbox}

\newcommand{\PP}{{\mathcal{P}}}

\newcommand{\eps}{{\epsilon}}

\newcommand{\Proj}{{{\mathtt{Proj}}^{\,}}}

\def\C{{\mathbb{C}}}

\def\PP{{\mathbb{P}}}
\def\Z{{\mathbb{Z}}}

\def\oo{{\mathcal O}}
\def\K{{\mathcal K}}

\def\D{{\mathcal D}}

\def\CO{{\mathcal O}}

\def\Gr{{\mathsf{Gr}}^{\mathsf{ad}}}

\def\ccirc{{}_{^\circ}}

\newcommand{\mmod}{\mathop{\op{\sf mod}}\nolimits}

\newcommand{\coh}{\mathop{\op{\sf coh}}\nolimits}
\newcommand{\qgr}{\mathop{\op{\sf qgr}}\nolimits}
\newcommand{\gr}{\mathop{\op{\sf gr}}\nolimits}
\newcommand{\tor}{\mathop{\op{\sf tor}}\nolimits}

\newcommand{\lra}{{\longrightarrow}}

\nc{\Fl}{\mathop{\rm Fl}\nolimits}

\nc{\newi}[1]{\noindent{\rm(#1)}}
\nc{\CC}{{\mathcal{C}}}
\nc{\CE}{{\mathcal{E}}}
\nc{\CF}{{\mathcal{F}}}
\nc{\CH}{{\mathcal{H}}}
\nc{\Tor}{\mathop{{\rm Tor}}\nolimits}
\nc{\CHom}{\mathop{\underline{\rm Hom}}\nolimits}
\nc{\CExt}{\mathop{\underline{\rm Ext}}\nolimits}
\renewcommand{\Hom}{\mathop{{\rm Hom}}\nolimits}
\nc{\Ext}{\mathop{{\rm Ext}}\nolimits}
\nc{\End}{\mathop{{\rm End}}\nolimits}
\nc{\Ker}{\mathop{\mathbf{Ker}}\nolimits}
\nc{\Rep}{\mathop{\mathbf{Rep}}\nolimits}
\nc{\Repm}{{\mathop{\mathbf{Rep}}\nolimits(\Gamma/\Gamma_\mu)}}

\nc{\Coker}{\mathop{\mathbf{Coker}}\nolimits}
\nc{\k}{\Bbbk}
\nc{\TV}{{\widetilde V}}
\nc{\CN}{{\mathcal N}}

\nc{\gtimes}{\otimes_\Gamma}
\newcommand{\ga}{{{\boldsymbol{\Gamma}}_{\!n}}}
\newcommand{\hh}{{\mathsf{H}}}

\nc{\G}{\Gamma}
\nc{\bplus}{\mbox{$\bigoplus$}}
\nc{\bbm}{{\mathbb{M}}}
\nc{\CW}{{\mathcal{W}}}
\nc{\gotimes}{{\otimes_{_\G}}}
\nc{\Symb}{{\mathop{\mathsf{Symb}}}}
\nc{\izp}{i^z_P}
\nc{\iwp}{i^w_P}
\nc{\tp}{{\tilde{p}}}
\nc{\tq}{{\tilde{q}}}

\newcommand{\Graff}{{\mathsf{Gr}}^{\mathsf{aff}}}
\newcommand{\Gradb}{{\mathsf{Gr}}^{\B}}
\newcommand{\Gradp}{{\mathsf{Gr}}^{\PPP}}
\newcommand{\Grn}{{\mathsf{Gr}^{\mathsf{ad}}_0}}

\newcommand{\Gradbp}{{\mathsf{Gr}}^{\B}_{p}}
\newcommand{\Gradpp}{{\mathsf{Gr}}^{\PPP}_{p}}
\newcommand{\Gradbq}{{\mathsf{Gr}}^{\B}_{q}}
\newcommand{\Gradpq}{{\mathsf{Gr}}^{\PPP}_{q}}

\newcommand{\MPPP}{{{\bf\mathcal{M}}_{\PPP}}}

\newcommand{\Gmu}{{\Gamma\cdot\mu}}

\begin{document}
\setlength{\parindent}{6mm}
\setlength{\parskip}{3pt plus 5pt minus 0pt}

\centerline{\Large {\bf  Wilson's Grassmannian and a Noncommutative
Quadric}}

\vskip 7mm
\centerline{\large {\sc {Vladimir Baranovsky, Victor Ginzburg and
Alexander
 Kuznetsov}}}

\vskip 7mm
\centerline{\it To Yuri Ivanovich Manin  on his 65-th birthday}
\bigskip
\begin{abstract}
{\footnotesize Let the group $\,\bmu\,$ of
$m$-th roots of unity act on the complex
line  by multiplication.
This gives a $\bmu$-action
on ${\mathsf {Diff}}$, the algebra of polynomial differential
operators on the line.
Following Crawley-Boevey and Holland, we introduce
a multiparameter  deformation, $\B_\tau$, of the smash-product
${\mathsf {Diff}}\#\bmu$. Our main result
provides  natural bijections between (roughly speaking) the following
spaces:
\vskip 1pt

(1)\; $\bmu$-equivariant version of
 Wilson's {\it adelic Grassmannian} of rank $r$;

(2)\; Rank $r$ projective $\B_\tau$-modules 
(with generic trivialization data);

(3)\; Rank $r$ torsion-free sheaves on a `noncommutative quadric' $\;\PPP$;

(4)\; Disjoint union of Nakajima quiver varieties for the cyclic
quiver with $m$ vertices.
\vskip 1pt

\noindent
The bijection between (1) and (2) is provided by a version
of Riemann-Hilbert correspondence between $\D$-modules and sheaves.
The bijections between (2), (3) and (4) were motivated by our previous work
\cite{BGK}. The resulting bijection  between (1) and (4)
reduces, in the very special case: $r=1$ and $\bmu=\{1\},$
to the partition of  (rank $1$) adelic Grassmannian
into a union of Calogero-Moser spaces, discovered by Wilson. This gives,
 in particular, a natural
and purely algebraic approach to  Wilson's result \cite{W}.}
\end{abstract}

{\centerline{\bf Table of Contents}
\vskip -5mm
$\hspace{20mm}$ {\footnotesize \parbox[t]{115mm}{\,

\hp${}_{}$\hp1.{ $\;\,\,$} {\tt Introduction}
\newline
\hp2.{ $\;\,\,$} {\tt Statement of Results}\newline
\hp3.{ $\;\,\,$} {\tt Kashiwara Theorem and De Rham Functor}\newline
\hp4.{ $\;\,\,$} {\tt $\B_\tau$-module Grassmanian and 
Sheaves on $\PPP$}\newline
\hp5.{ $\;\,\,$} {\tt Monads and Quiver Varieties}\newline
\hp6.{ $\;\,\,$} {\tt Projective $\B_\tau$-modules}\newline
\hp7.{ $\;\,\,$} {\tt Appendix A: Formalism
of Polygraded Algebras}\newline
\hp8.{ $\;\,\,$} {\tt Appendix B: The Geometry of $\PPP$}\newline
}}}

\section{Introduction}
\setcounter{equation}{0} 
Nakajima \textit{quiver varieties} can be viewed, according to
our previous paper \cite{BGK}, as spaces parametrizing torsion-free
sheaves on a ``noncommutative plane". In the simplest case this yields a
relation, first observed by Berest-Wilson~\cite{BW}, between Calogero-Moser 
spaces and projective modules over the first
Weyl algebra, $\D(\C)$, of polynomial differential
operators on the line $\C$. The approach to this result (and to its
`quiver generalizations') used in \cite{BGK} was purely algebraic
and totally different from the approach in~\cite{BW}. The latter involved
a non-algebraic {\it Baker function} and was based heavily on the earlier
remarkable discovery by Wilson \cite{W} of a connection between
an {\it adelic Grassmannian} and  Calogero-Moser
spaces.

In this paper we reverse the
logic used by  Berest-Wilson
 and
 {\it explain} (rather than {\it exploit}) the 
 connection between
adelic Grassmannians and  Quiver varieties
by means of `noncommutative algebraic geometry',
using the strategy of
 \cite{BGK}.
\footnote{During the preparation of the present paper (which was first supposed
to be part of  \cite{BGK}) another paper by Berest-Wilson
appeared, see \cite{BW2}. Our  approach is similar to that of \cite{BW2}
(we treat more general
case of `higher rank' and $\bmu$-equivariant
projective modules). However, even in the rank 1 case,
 in \cite{BW2} the authors do 
not provide an independent proof of the bijection between Calogero-Moser 
spaces and projective modules; instead they construct a map inverse
to the map constructed in \cite{BW} assuming the latter is already known
to be a bijection. An independent direct proof  of the bijectivity  in
the rank 1 case
was
obtained in the Appendix to  \cite{BW2}  by M. Van den Bergh, who
 used some results of  \cite{BGK}. 

We emphasize that, for the reasons
  explained at the end of the
Introduction below, it seems to be impossible to
extend the approach
of \cite{BW2} (connecting  the adelic
Grassmannian with rank 1 sheaves on a non-commutative
${\mathbb{P}}^2_\tau$) to the higher rank case without replacing
${\mathbb{P}}^2_\tau$ by a non-commutative  surface which fibers
over ${\mathbb{P}}^1$,  like the  noncommutative quadric
$\PPP$ that we are using in the present
paper.}

 Our first key observation is
that each point of adelic Grassmannian can be viewed
as a `constructible sheaf' on the line built up from
`infinite-rank' local systems. This way, the correspondence between
projective (not {\it holonomic}\,!) $\D(\C)$-modules
 and points of the  adelic Grassmannian becomes
nothing but  (a non-holonomic version of)
the standard De Rham functor between $\D$-modules and
constructible sheaves on the line\footnote{More generally,
our construction
of De Rham functor yields a similar correspondence between projective
$\D$-modules on any smooth algebraic curve  $X$
and points of an appropriately defined adelic Grassmannian 
attached to the curve (in that case a non-commutative version of
projective completion of $T^*X$ should play the role of  $\PPP$).
The case of an elliptic curve
seems to be especially interesting; we hope to discuss it elsewhere.}.

The  De Rham correspondence works equally well in a  more general context
of equivariant $\D$-modules with respect to a natural
action on the line $\C$ of the group $\bmu$ of $m$-th roots of unity,
by multiplication.
Giving a $\bmu$-equivariant $\D$-module is clearly the
same thing as giving a module over $\D(\C)\#\bmu,$ the smash-product of
 $\D(\C)$  with the group  $\bmu$,
acting on  $\D(\C)$ by algebra automorphisms. Note that
in \cite{BGK} {\it any}, not only cyclic, finite
group $\Gamma \subset SL_2(\C)$ of automorphisms
of the Weyl algebra was considered. In order
to have a   De Rham functor, however,
one needs to specify a
 standard holonomic $\D$-module of
`regular functions'. The choice of such a $\D$-module
breaks the $SL_2(\C)$-symmetry of the 2-plane  formed
by the generators of the first Weyl algebra.
Thus, the group  $\Gamma \subset SL_2(\C)$ has
to have an invariant line in $\C^2$. This leaves us with the only choice
$\Gamma=\bmu$.

Below, we will be working
not only with the  algebra $\D(\C)\#\bmu$,
but with a multi-parameter deformation 
\begin{equation}\label{D}
\B_\tau = \C\langle x,y\rangle\#\bmu^{\,} \big/^{\,}
 \langle [y,x]=\tau \rangle.
\end{equation}
 of
that algebra introduced by Crawley-Boevey 
and Holland. Here $\tau$ (= deformation parameter) is  an
arbitrary element 
in the  group algebra $\C[\bmu]$, and $\C \langle x, y \rangle$ 
stands for the free $\C$-algebra
of noncommutative polynomials in two
variables $x, y$.
Once  a De Rham  functor  between
projective $\B_\tau$-modules  and points of an
 adelic Grassmannian is established one can
construct
a `Wilson type' connection between the
adelic Grassmannian and  Quiver varieties as follows.
First,  view a projective $\B_\tau$-module
as a vector bundle on an appropriate noncommutative
plane ${\mathbb{A}}^2_\tau$. Next,  extend (see  \cite{BGK})
this vector bundle to a (framed) torsion-free sheaf on a {\it completion}
$X_\tau\supset {\mathbb{A}}^2_\tau$,
a `noncommutative
projective surface'. Finally, we use a 
description of framed torsion-free sheaves on $X_\tau$
in terms of {\it monads} (i.e. in terms of linear algebra data)
developed in \cite{BGK} to obtain a parametrisation
of  projective $\B_\tau$-modules by points of 
certain Quiver varieties.

There are several possible choices for a `compactification' $X_\tau$
of the  noncommutative
plane~${\mathbb{A}}^2_\tau$. In \cite{BGK} we  used $X_\tau=\PP^2_\tau,$
a  noncommutative version of projective plane.
In the present paper we choose 
another `compactification' of  ${\mathbb{A}}^2_\tau$,
a  noncommutative version, $\PPP$, of two-dimensional quadric.
This choice is essential for our present purposes. Our construction
of the extension of a $\B_\tau$-module to a torsion free sheaf
on $\PPP$ doesn't behave well enough in the case of $\PP^2_\tau$.
On the other hand, the relation of sheaves on $\PPP$ to quiver 
varieties is {\it a posteriori} equivalent to
the one used in \cite{BGK}, since the two noncommutative spaces
$\PPP$ and  $\PP^2_\tau$ can be obtained from each other
by "blowing up" and "blowing down" constructions.
We will  indicate the idea of such a construction
in a Remark (above Theorem \ref{betaiso}) and it
will be hinted  there how a
bijection between torshion-free sheaves on 
 $\PPP$ and  on $\PP^2_\tau$ can be  established via a  
noncommutative version of Fourier-Mukai transform,
see~\eqref{radon}.

Our  results generalize (and, hopefully, put in context) the results of
Wilson \cite{W} in two ways. First, we incorporate a $\bmu$-action.
Second, Wilson only considered the rank 1 case,
that is the case of rank 1 sheaves on $\PPP$ in the terminology
of the present paper. 
In Wilson's situation the 
\textit{whole} adelic Grassmanian gets partitioned
into a disjoint union of Calogero Moser spaces. 
In the more general setup of arbitrary rank $r\geq 1$
this is no longer true for two reasons. First, in our
definition of the adelic Grassmanian we drop the  "index zero"
condition of
$\bigl(\cite{W}, 2.1(ii)\bigr)\,$ (it has to do with replacing the group
$SL_r$ by $GL_r$).
 This makes our version of
 adelic Grassmanian  somewhat larger than 
that of \cite{W}. The geometric consequence of  "index zero" condition 
is (in our language)
that the restriction of a  coherent sheaf to the line 
$\PP^1\times\{\infty\} \subset \PPP$
is a vector bundle with the vanishing first Chern class: $c_1=0$. 
In the rank~1 case considered by Wilson, 
any such bundle is necessarily
trivial, while this is clearly not true for higher ranks.
Thus, our main result says that, for any $m\geq 1$ and  $r\geq 1$,
 the part of ($\bmu$-equivariant)
rank $r$ adelic Grassmanian formed by sheaves {\it trivial} on
the line
$\PP^1\times\{\infty\} \subset \PPP$ can be partitioned into
 a disjoint union of quiver
varieties of type $\mathbf{A_{m-1}}$.
\medskip

\begin{remark} We would like to end this Introduction by bringing
attention
of the reader to 
a surprising correspondence resulting from comparing  
 \cite[Theorem 1.3.12]{BGK} with \cite[Theorem 1.13]{EG}.
Specifically, let $\Gamma\subset SL_2(\C)$ be a finite subgroup,
and $\B_\tau(\Gamma)$ the corresponding coordinate ring of
the corresponding $\Gamma$-equivariant `non-commutative plane'.
Recall that according to \cite{BGK},
there is a  partition of
the moduli space of projective $\B_\tau(\Gamma)$-modules $N$ such that 
$[N]=R$  into a disjoint union,
according to the `second Chern class' $c_2(N):=n\in \Z$.
On the other hand, given $n\ge 1$, let
$\,\ga:=S_n\ltimes (\Gamma\times \Gamma\times\ldots\times \Gamma)\subset
Sp(\C^{2n})$
denote the corresponding wreath-product and let
$\hh_{0,\tau}(\ga)$ be the {\it Symplectic Reflection algebra}
attached to $\ga$, see \cite[p.249]{EG}.
\footnote{In \cite{BGK} we use the notation  ${\mathcal{B}}_\tau$
instead of
$\B_\tau(\Gamma)$, and in \cite{EG} we write `$c$' for what we denote by
$\tau$ in \cite{BGK} (and in the present paper).}  
Further,
Theorem 1.13 of \cite{EG} shows that representation theories of the
algebras $\B_\tau(\Gamma)$ and $\hh_{0,\tau}(\ga)$ are related by the following
{\bf mysterious bijection}:

{\small
$$
\left\lbrace
\begin{array}{l}\mbox{Isomorphism classes of finitely generated 
{\bf projective}}\\
\mbox{$\B_\tau(\Gamma)$-modules $N$ such that 
$[N]={\mathsf{triv}}$  and $c_2(N)=n$}
\end{array}
\right\rbrace
\;\;\simeq\;\;
\left\lbrace
\begin{array}{l}\mbox{Isomorphism classes of}\\
\mbox{{\bf simple} $\hh_{0,\tau}(\ga)$-modules}
\end{array}
\right\rbrace
$$ 
}
Finding a direct conceptual explanation of the bijection above presents
a challenging problem. We remark that even in the case of the
trivial group $\Gamma$, where the moduli space on each side 
reduces, as a variety, to the Calogero-Moser space,
the bijection is still completely unexplained.
\end{remark}

\bigskip
\noindent
{\bf Acknowledgements.} {\footnotesize
We are  indebted to Sasha Beilinson for
some very useful remarks.
The third author was partially supported by RFFI grants
99-01-01144 and 99-01-01204. Also, he would like 
to express his gratitude to the University of Chicago, 
where the major part of this paper was written.}

\section{Statement of Results}
\setcounter{equation}{0}
   From now on, let $\Gamma = \bmu$ denote the group of $m$-th roots of
unity, and $\CG$  its group algebra. We fix an embedding
$\Gamma = \bmu\into SL_2(\C)$, and let
$L$ denote the tautological two-dimensional
representation of $\Gamma$
arising from the embedding. We have: $L\cong\eps\oplus\epss$,
where $\eps$ is a primitive
character of $\Gamma$. Let $\{x,y\}$ be a basis of $L^*$ compatible
with the above direct
sum decomposition, such that  $\Gamma$ acts on $x$ by $\eps$ and on $y$ by
$\epss$. Write $\C[x]$ for the polynomial algebra on the line
with coordinate $x$,
and $\C(x)$ for the corresponding field of rational functions.
We form the smash-product algebras
$$
\CG[x] \,:=\, \C[x]\#\Gamma\quad\text{and}\quad
\CG(x) \,:=\, \C(x)\#\Gamma.
$$
The standard embedding $\CG\into \C[x]\#\Gamma$ 
makes $\CG[x]$ into a $\Gamma$-bimodule
via left and right multiplication by $\G$.
There is a similar $\Gamma$-bimodule structure on $\CG(x)$.

Choose and fix a finite-dimensional $\Gamma$-module $W$.
There is a natural $\C(x)^\G$-action on
 $W\gotimes\CG(x)$
 given by $p:\, w\otimes f \longmapsto p\cd(w\otimes f)
:=w\otimes (p\cd f).$ 

\begin{definition}
 A $\Gamma$-invariant vector subspace $\CU\subset W\gotimes\CG(x)$ is
called {\sl {primary}}\linebreak
 {\sl {decomposable}}\footnote{this notion is due
to Cannings-Holland \cite{CaH}.}
 if the following two conditions
hold

(a) There exists a $\Gamma$-invariant polynomial $p=p(x)$ such that
$$
p\cd (W\gotimes \CG[x])\;\subset\;\CU\;\subset\; 
\mbox{$\frac{1}{p}$}\cd(W\gotimes \CG[x]).
$$

(b) If $p(x) = \prod_{\mu}(x-\mu)^{k_\mu}$ then the subspace
on the left below is compatible with the direct sum decomposition
on the right (i.e. LHS = sum of its intersections with the
direct summands on the RHS):
$$
\frac{\CU}{p\cd (W\gotimes\CG[x])}\; \subset\;
\frac{\po(W\gotimes \CG[x])}{p\cd (W\gotimes\CG[x])}\; =\;
\bigoplus_{\mu}\;
\frac{(x-\mu)^{-k_\mu}W\gotimes\CG[x]}{(x-\mu)^{k_\mu}W\gotimes\CG[x]}\;\;.
$$

Define an \emph{adelic Grassmanian} $\Gr(W)$ to be the set of all
primary decomposable subspaces $\CU\subset W\gotimes\CG(x)$.

\end{definition}

Our first goal is to relate  the adelic Grassmanian to
 modules
over some 
noncommutative algebra. To that end, we fix an  element
$\tau\in \CG$ and consider the algebra $\B_\tau$, see \eqref{D},
to be denoted $\B$ in the future.
Let  $\Brat$ be the localization of the algebra
$\B$ with respect to the
multiplicative system $\,{\C[x]^\G\smallsetminus\{0\}}\,$
of all nonzero $\Gamma$-invariant polynomials in $x$. This localization has
a natural  algebra structure extending that on $\B$.
Note further that we have a natural embedding
$\CG(x)\into \Brat$ and, moreover,
this embedding yields a vector space isomorphism:
$\CG(x)\iso \Brat/\Brat\cd y.$
We make $\CG(x)$ into a left $\Brat$-module by
 transporting the obvious 
 $\Brat$-module structure
on $\Brat/\Brat\cd y$ via the bijection above.
For $\tau=1$, this
 reduces essentially to the standard
action on  $\CG(x)$ by differential operators.

For any  $\G$-module $W$,
the space  $W\gotimes\Brat$ has an obvious structure
of a projective right  $\Brat$-module.

\begin{definition}\label{group}
Let $\Gw=GL_{\Brat}(W\gotimes\Brat)$ be the group
of all (invertible) right $\Brat$-linear automorphisms
of the  $\Brat$-module $W\gotimes\Brat$.
\end{definition}

We define a
{\it left} $\Gw$-action on the vector space  $W\gotimes\CG(x)$
and on the adelic Grassmanian $\Gr(W)$
as follows. 
First, observe that the natural  left
$GL_{\Brat}(W\gotimes\Brat)$-action on
 $W\gotimes\Brat$ commutes with right multiplication
by `$y$', therefore keeps the subspace
${W\gotimes\Brat\cd y}\subset W\gotimes\Brat$ stable.
Hence, there is a well-defined  left $\Gw$-action on
the quotient $(W\gotimes\Brat)\big/(W\gotimes\Brat\cd y)$.
Further, since left $\Gamma$-action commutes with right 
$\Brat$-action, and since $W$ is a projective $\Gamma$-module,
it follows that we have an isomorphism
$$
\psi:\; W\gotimes\CG(x) \iso (W\gotimes\Brat)\big/(W\gotimes\Brat\cd y),
$$
induced by the isomorphism $\CG(x)\iso \Brat/\Brat\cdot y$
considered two paragraphs above. We define the left
$\Gw$-action on $W\gotimes\CG(x)$ by transporting the left
$\Gw$-action on the quotient $(W\gotimes\Brat)\big/(W\gotimes\Brat\cd y)$
via the  isomorphism $\psi$.
It is straightforward to verify that
elements of $\Gw$ take primary decomposable subspaces
of $W\gotimes\CG(x)$
into primary decomposable subspaces. This gives a
canonical left $\Gw$-action on 
the adelic Grassmanian~$\Gr(W)$.

To formulate our first result, recall that
there is a canonical isomorphism (due to Quillen) 
of Grothendieck
 K-groups: $ K(\CG)\iso K(\B)$ induced by the
 functor:
$W \mapsto W\gotimes\B$.
We write $[N]\in K(\CG)$ for the image of the class of
a $\B$-module $N$ under the inverse isomorphism
$K(\B)\to K(\CG)$, and let $\dim: K(\CG)\to\Z$
denote the dimension homomorphism.

Further,
there is a distinguished finite collection
of codimension one hyperplanes
in the vector space $\CG$, called the {\it root} hyperplanes.
One way to define these  hyperplanes is to use
  McKay correspondence. The latter associates to the
cyclic group $\G=\bmu\subset SL_2(\C)$ an affine Dynkin graph
of type $\widetilde{\mathbf{A_{m-1}}}$ such that
the underlying vector space of the group algebra
 $\CG$ gets identified  with the dual of the
$\C$-vector space generated by simple roots of the corresponding
affine root system. In particular, every root gives a hyperplane
in the vector space $\CG$.

\begin{definition}\label{def_gen}
An  element $\tau\in \CG$ is called \emph{generic} if
it does not belong to  any  root hyperplane in $\CG$.
\end{definition}

Our first theorem below is a noncommutative analogue of a well-known
result due to A.~Weil, providing a description of algebraic vector
bundles on an algebraic curve in terms of an adelic double-coset
construction.

\begin{theorem}\label{prb}
Assume that $\tau$ is generic.
Let $W$ be a $\Gamma$-module with $\dim W=r$.
The set of isomorphism classes of projective  (right) $\B$-modules $N$
such that $\dim[N]=r$ is in a canonical bijection with
the coset space $\Gw\backslash\Gr(W)$.
\end{theorem}

To explain the main ideas involved in the proof of Theorem~\ref{prb}
we need the following

\begin{definition}
A right $\B$-submodule $N\subset W\gotimes\Brat$ is called \emph{fat}
if there exists a $\Gamma$-invariant polynomial $p(x)$ such that
$p\cd (W\gotimes \B)\subset N \subset \frac{1}{p}\cd(W\gotimes \B)$.

Let $\Gradb(W)$ be the set of all
fat right $\B$-submo\-dules $N \subset W\gotimes\Brat$.
\end{definition}

Now, the proof goes as follows. First, we check that for 
generic~$\tau$ any projective right $\B$-module can be 
embedded into $W\gotimes\Brat$ as a fat $\B$-submodule.
The embedding is unique up to a $\Gw$-action. Then, it 
remains to relate the Grassmanians $\Gradb(W)$ and $\Gr(W)$.
To this end,
recall that we have equipped the space $\CG[x]$ with a canonical 
structure of left $\B$-module, that clearly commutes with right 
$\G$-action by multiplication. 

We introduce a De Rham functor $\DR$ from
the category of {\it right} $\B$-modules to the category of right
$\Gamma$-modules as follows
$$
N\longmapsto \DR(N) := N\otimes_\B\CG[x].
$$
Given a non-zero polynomial $p\in \C[x]^\G,$ 
we write $p\cd (W\gotimes \B) := \bigl(W\gotimes (p\cd\B)\bigr).$
The space 
$p\cd (W\gotimes \B) $
has an obvious right $\B$-module structure, and we have:
$$
\DR\bigl(p\cd (W\gotimes \B)\bigr) =
p\cd (W\gotimes \B)\otimes_\B\CG[x] =
p\cd (W\gotimes\CG[x])\;.
$$
Hence the De Rham functor takes any fat $\B$-submodule of
$W\gotimes\Brat$ to a $\Gamma$-invariant vector subspace
$\CU\subset W\gotimes\CG(x)$, such that
$p\cd (W\gotimes\CG[x])\subset\CU\subset \po(W\gotimes \CG[x])$.
Moreover, one can check that this space is primary decomposable.

\noindent
{\bf Example.\;}
Assume that $\Gamma=\{1\}$ is trivial, i.e. that $m=1$.
Then $\CG=\C$ and $\tau\in \CG=\C$ 
is generic if and only if  $\tau\ne0$. In
this
case the algebra $\B$ is isomorphic to $\D(\C)$,
the algebra of differential operators
on  the line $\C$, and the algebra $\Brat$ is isomorphic to the algebra of
differential operators
with rational coefficients. The functor $\DR$ becomes
 the standard De Rham functor.
\bigskip

We have seen that the De Rham functor yields a map
$\DR: \Gradb(W)\to\Gr(W)$.
We also define a map: $\Gr(W)\to\Gradb(W)$ (in the opposite direction)
as follows.
Let $\CU\subset W\gotimes\CG(x)$ be a primary decomposable
subspace. Then the left $\Brat$-action map $\Brat\otimes\CG(x)\to\CG(x)$
induces, after tensoring with $W$ and restricting to $\CG[x]$, a linear map
$a: W\gotimes\Brat\otimes_\C\CG[x]\too
W\gotimes\CG(x)$. Let $\Diff(\CU)$ denote the set of all elements $u\in
W\gotimes\Brat$,  such that $a(u\otimes\CG[x])\subset \CU$.
It is easy to show that $\Diff(\CU)\subset W\gotimes\Brat$ 
is a fat $\B$-submodule. Thus, we obtain a map $\Diff:\Gr(W)\to\Gradb(W)$.

\begin{theorem}\label{DRiso}
For generic $\tau$, the maps $\DR$ and $\Diff$
give mutually inverse bijections
$\Gradb(W)\cong\Gr(W)$.
\end{theorem}

Our next step is to interpret
the space $\Gradb(W)$ using the formalism of
noncommutative geometry. To this end, we introduce
 the algebra
\begin{equation}\label{qu}
\Qu = \C\langle x,z,y,w\rangle\#\Gamma \Big/
\Big\langle [x,z] = [y,z] = [z,w] = [y,w] = [x,w] = 0,\
[y,x]=\tau\cd zw \Big\rangle\;,
\end{equation}
where for any $\gamma\in\G$ we put:
\begin{equation}\label{gxyzw} 
\gamma\cdot x\cdot\gamma^{-1}=\epsilon(\gamma)\cd x,\quad
\gamma\cdot y\cdot\gamma^{-1}=\epsilon^{-1}(\gamma)\cd y,\quad
\gamma\cdot z\cdot\gamma^{-1}=z,\quad
\gamma\cdot w\cdot\gamma^{-1}=w.
\end{equation}
Define a bigrading: $\Qu=\oplus_{i,j\ge 0} \Qu_{i,j}$ on
the algebra $\Qu$ by letting
$\deg x = \deg z = (1,0)$, $\deg y = \deg w = (0,1),$
and $\deg\gamma = (0,0)$ for any $\gamma\in\Gamma$.
Thus,
$\Qu_{0,0} = \CG$.

When $\Gamma$ is trivial and  $\tau = 0$ the algebra $\Qu$ 
reduces to the standard
bigraded
algebra associated to the quadric $\PP^1 \times \PP^1$ and a pair of
line bundles $L_1 = \oo (1, 0)$;  $L_2 = \oo(0, 1)$, 
that is for any $i,j\geq 0$, we have $\Qu_{i, j} =
H^0(\PP^1 \times
\PP^1, L_1^{\otimes i} \otimes L_2^{\otimes j})$. In this case,
 the category of coherent
sheaves of
$\PP^1 \times \PP^1$ can be described as a quotient category of the category
of
bigraded $\Qu$-modules (see \S4).

In the general case of a nontrivial group
$\Gamma$ and arbitrary $\tau$
 a similar quotient category
construction may  still  be applied formally to
the  bigraded ring $\Qu$. 
Following  \cite{AZ}, see also \cite{BGK} and \S\S7-8 below,
we will  view objects of the resulting  quotient category  as 
coherent sheaves on a
``noncommutative
quadric'' $\PPP$.

Note that, in the commutative case $\tau=0$, 
the equations $z =0$ and $w = 0$ give rise to
two embeddings
$i_z: \PP^1_z \hookrightarrow \PP^1 \times \PP^1$ and $i_w: \PP^1_w
\hookrightarrow
\PP^1 \times \PP^1$
of the corresponding factors $\PP^1$.
Thus, one may
 consider the restriction
functor $i^*_z$ taking
coherent sheaves on $\PP^1 \times \PP^1$ to coherent sheaves on
$\PP^1_z$
(and also  consider coherent sheaves on $\PP^1 \times \PP^1$
which are trivialized
in some formal neighbourhood of $\PP^1_z$).
In \S4 and Appendix A, we
show how to extend all relevant concepts to the noncommutative case. The
homogeneous
coordinate rings of $\PP^1_z$ and $\PP^1_w$ will be replaced by $\C[y,
w]\#\Gamma$ and
$\C[x, z] \#\Gamma$, respectively. 
The  latter 
algebras are only slightly noncommutative in the sense
that the corresponding quotient
categories of  `noncommutative coherent 
sheaves' are nothing but the categories of
$\Gamma$-equivariant  coherent  sheaves on $\PP^1$,
the  ordinary 
(commutative) projective line.
 This leads to the following provisional
(see Definitions \ref{trivdef1} and \ref{trivdef2} for details)
\begin{definition}\label{nc-pre} 
Let $\Gradp(W)$ be the set of all {\rm(}equivalence classes of\/{\rm)} 
coherent sheaves $E$ on $\PPP$ trivialized in a formal 
neighbourhood of $\PP^1_z$ and such that 
$i^*_z E \simeq W \otimes_{\Gamma} \oo_{\PP^1_z}$.
\end{definition}

Note that in the commutative situation we have:
$\PP^1 \times \PP^1 \smallsetminus
(\PP^1_z \cup \PP^1_w)= {\mathbb A}^2$
is an affine plane with coordinates $x, y$. In the noncommutative case we
have an  algebra isomorphism
$$
\B \simeq \Qu\big/\bigl(( z- 1) \Qu + (w - 1) \Qu\bigr)\;\;.
$$
Therefore, the algebra $\B$ may be viewed as coordinate ring
of a noncommutative affine plane $j:{\mathbb A}^2_\tau\to
\PPP$. This gives rise to a restriction functor
$j^*:\coh(\PPP)\to\mmod(\B)$ taking the
category of coherent sheaves on $\PPP$ to $\B$-modules. It is easy to see
that $j^*$ takes
any sheaf trivialized in a neighborhood of $\PP^1_z$ to
a fat $\B$-submodule of $W\gotimes\Brat$.

\begin{theorem}\label{jstariso}
The `restriction' functor $j^*$ induces a bijection $\Gradp(W)\iso
\Gradb(W)$.
\end{theorem}
We write $j_{!*}$ for an inverse of the bijection $j^*$.
\bigskip

The fourth (and the last) infinite Grassmanian considered
in this paper is an \emph{affine Grassmanian} $\Graff(W)$ introduced
below. 

\begin{definition}
A right $\CG[x]$-submodule $\CW\subset W\gotimes\CG(x)$ is called \emph{fat}
if there exists a $\Gamma$-invariant polynomial $p(x)$ such that
$p\cd (W\gotimes\CG[x])\subset\CW\subset \po(W\gotimes \CG[x])$.

Define $\Graff(W)$ to be the set of all
fat $\CG[x]$-submodules in $W\gotimes\CG(x)$.
\end{definition}

\noindent
It is clear that a (right) $\G$-stable subspace $\CW\subset W\gotimes\CG(x)$
is a fat $\CG[x]$-submodule if and only if $\CW$ is a finitely-generated
$\C[x]$-submodule such that $\CW\otimes_{_{\C[x]}}\,\C(x)=W\gotimes\CG(x)$.
Thus, a fat $\CG[x]$-submodule may be viewed
as a $\G$-stable {\it lattice} in the $\C(x)$-vector space
$W\gotimes\CG(x)$. For this reason we refer to $\Graff(W)$
 as the
\emph{affine Grassmanian}. The standard relation between
loop-Grassmannians and vector bundles on the Riemann sphere, see e.g. \cite{PS},
shows that the space
$\Graff(W)$ can also
be interpreted as the set of all $\Gamma$-equivariant
vector bundles
on $\PP^1$ trivialized in a Zariski neighbourhood of 
the point $\infty \in \PP^1$
(cf. Definition~\ref{trivdef1}),
with $W$ being
a fiber at $\infty$. Thus the  restriction functor $i^*_w$ takes
$\Gradp(W)$ to
$\Graff(W)$. 

We observe  that the group $\Gw$, see Definition \ref{group},
acts naturally on each of the  Grassmannians:
$\Gr(W)\,,\,\Gradb(W)$, and $\Gradp(W)$.
Specifically, the action on $\Gr(W)$ has been defined earlier,
and the  action on $\Gradb(W)$
is induced by
the corresponding $\Gw$-action on $W\gotimes\Brat$.
The action on $\Gradp(W)$
arises from  $\Gw$-action on the direct
system formed by the sheaves $W\gotimes\CO(n,0)\,,\,n=1,2,\ldots$.
Observe further that the affine Grassmanian $\Graff(W)$ has an action
of the subgroup $GL_{\CG(x)}(W\gotimes\CG(x))\subset\Gw$
(the group $\Gw$ itself does not act on $\Graff$ since
it does not preserve the condition to be a $\CG[x]$-submodule).

All the objects and the maps we have
introduced so far are incorporated in
the following  diagram
\begin{equation}\label{gwequi}
\xymatrix{
\Gr(W)     \;     \ar@<.5ex>[rr]^{\Diff}  &&
\;\Gradb(W)  \;     \ar@<.5ex>[ll]^{\DR} \ar@<.5ex>[rr]^{j_{!*}} 
&&
\;\Gradp(W) \;      \ar@<.5ex>[ll]^{j^*} \ar[r]^<>(.5){i_w^*}& 
\; \Graff(W)\;\;,
}
\end{equation}
\medskip

\noindent
{\bf {Principal symbol map `$\Symb$':}}\;\;
The algebra $\B$ has a natural  increasing
 filration: $\CG[x]=F_0\B \subset F_1\B\subset F_2\B\subset\ldots,$
where $F_k\B$ is the $\CG[x]$-submodule generated
by $\{1,y,\ldots,y^k\}.$ This filtration by
the `order of differential operator'
extends canonically to a similar   filtration
$\CG(x)=F_0\Brat\subset F_1\Brat\subset F_2\Brat
\subset\ldots,$ on the algebra $\Brat$,
and for the corresponding associated graded algebras we have:
$\gr^F\B\simeq\CG[x,y]$ and $\gr^F\Brat\simeq\CG(x)[y]$.
The filtration on  $\Brat$  also induces  an  increasing filtration,
$F_k(W\gotimes\Brat):= W\gotimes F_k\Brat,$ on the
$\Brat$-module $W\gotimes\Brat$ such that
$\gr^F(W\gotimes\Brat)\;\simeq\;\bigl(W\gotimes\CG(x)\bigr)[y].$

Now, given  a $\B$-submodule $N\subset W\gotimes\Brat$,
we put
\begin{align*}
\Symb(N):= \{ f\in W\gotimes\CG(x)\enspace \big|&\enspace
f\cdot y^k + u_{k-1} \in N,\\ 
&\text{ for some $k\in\Z$ and some
$u_{k-1}\in F_{k-1}(W\gotimes\Brat)$}\}\;.
\end{align*}
This is a $\CG(x)$-submodule in $W\gotimes\CG(x)$ that can be
equivalently defined as follows.

Right multiplication by `$y$' gives rise
to a direct system of bijective maps: \linebreak
$\CG(x)\iso\CG(x)\cd y\iso\CG(x)\cd y^2\iso\ldots,$
and it is clear that this yields isomorphisms
\begin{align}\label{dir}
&\CG(x)=\CG(x)\cd y^0\enspace\iso\enspace\dlimk\;\CG(x)\cd y^k
\quad\text{and}\\
&W\gotimes\CG(x)=W\gotimes\CG(x)\cd y^0
\enspace\iso\enspace\dlimk\;\bigl(W\gotimes\CG(x)\bigr)\cd y^k\;.\nonumber
\end{align}
Let $\gr^FN \subset \gr^F(W\gotimes\Brat)
\simeq W\gotimes\CG(x)[y]$ denote the associated graded of $N$
with respect to the induced filtration $F_\bullet N:= N\cap 
F_\bullet(W\gotimes\Brat)$,
and form the direct system:
$\gr^F_0N\to\gr^F_1N\to\gr^F_2N\to\ldots,$
induced by the $y$-action on $\gr^FN$ (which is not necessarily bijective).
Using the identification provided by \eqref{dir},
one has:
 $\Symb(N) =\dlim\;\gr^F_kN$.
It is clear that  the RHS is a $\CG[x]$-submodule in $\CG(x)$.
The assignment: $N\mapsto \Symb(N)$ gives a (discontinous)
map $\Symb:\Gradb(W)\to\Graff(W)$.
 
Let $\sigma$ denote the composite map, see \eqref{gwequi}
$\dis\sigma:\,\Gr(W)  \stackrel \Diff\too\Gradb(W)
\stackrel {j_{!*}}\too\Gradp(W)$ 
$\stackrel{i_w^*}\too\Graff(W).\,$ The 
following  is an enriched version of diagram
\eqref{gwequi}
\begin{theorem} \label{symb-restr}
Assume that $\tau$ is generic. Then we have a
 commutative diagram
\begin{equation}\label{pic}
\xymatrix{
\Gr(W)          \ar@<.5ex>[rr]^{\Diff} \ar[rrd]_{\sigma} &&
\Gradb(W)       \ar@<.5ex>[ll]^{\DR} \ar@<.5ex>[rr]^{j_{!*}} \ar[d]^{\Symb}
&&
\Gradp(W)       \ar@<.5ex>[ll]^{j^*} \ar@<.5ex>[lld]^{i_w^*} \\
&& \Graff(W)
}
\end{equation}
where the maps
$\DR\,,\,\Diff\,,\,j_{!*},$ and $j^*$ are
 $\Gw$-equivariant bijections, and 
the maps $i_w^*\,,\, \Symb,$ and $\sigma$ are
$GL_{\CG(x)}(W\gotimes\CG(x))$-equivariant.
\end{theorem}
\medskip

Finally, we explain how Quiver varieties enter the picture.
Given a 
pair of finite dimensional
 $\Gamma$-modules $W$, $V,$ and an element $\tau \in
\CG$ define, following Nakajima, 
a locally closed subvariety of \textit{quiver data}:
\begin{equation}\label{bbm}
\bbm^\tau_{_\G}(V,W) \;:=\;
\big\{(B,I,J)\in \Hom_{_\G}(V,V\otimes_\C L)
\;\bplus\;\Hom_{_\G}(W,V)\;\bplus\;\Hom_{_\G}(V,W)\big\}
\end{equation}
formed by the triples $(B,I,J)$ satisfying the following two
conditions:
\begin{description}
\item[Moment Map Equation:] $[B,B] + IJ = \tau\big|_V$;
\item[Stability Condition:] if $V'\subset V$ is a $\Gamma$-submodule
such that $B(V')\subset V'\otimes L$ and
\newline $I(W)\subset V'$ then $V'=V$.
\end{description}

\noindent
The group $G_{_\G}(V)=GL(V)^\G$ of $\Gamma$-equivariant automorphisms
of $V$ acts on $\bbm^\tau_{_\G}(V,W)$ by the formula:
$
g(B,I,J) = (gBg^{-1},gI,Jg^{-1}).
$
Note that this 
 $G_{_\G}(V)$-action is free, due to the stability condition.

\begin{definition}\label{quiverdef} Let
$\FM_{_\G}^\tau(V,W)=\bbm^\tau_{_\G}(V,W)/G_{_\G}(V)$
be the Geometric Invariant Theory quotient, called a
(Nakajima) quiver variety.
\end{definition}

The affine Grassmanian $\Graff(W)$ has a marked point $\CW_0 = W
\otimes_{\Gamma} \C\Gamma[x]$ 
corresponding to the $\Gamma$-equivariant sheaf $W \otimes_{\Gamma}
\oo_{\PP^1_w}$ with
its natural trivilization in the Zariski neighbourhood of infinity.

\begin{theorem}\label{framedsheaves}
Let $\tau$ be  generic. Then there is a canonical bijection
$$\bigsqcup\nolimits_V\;\FM_{_\G}^\tau(V,W)\;\cong\;\sigma^{-1}(\CW_0)
\;\subset\;\Gr(W)\;\;,
$$
where $V$ runs through the
set of isomorphism classes of all finite dimensional
$\Gamma$-modules.
\end{theorem}

\noindent
{\bf Example.\;} Let the group $\Gamma$ be trivial and  $\tau\ne 0$.
If $W=\C$ and $V=\C^n$, 
then the corresponding quiver variety is isomorphic to 
the Calogero-Moser space $CM_n$. Further,
 the affine Grassmanian
reduces to the coset space $\C(x)^\times/\C^\times$.
Moreover, the subset
$\Grn(\C):=\sigma^{-1}(\CW_0)\subset\Gr(\C)$ 
consists of primary decomposable subspaces of `index
zero' (in the sense of \cite[2(ii)]{W}) at every point. Thus,  Theorem
\ref{framedsheaves} implies in this case 
 Wilson's Theorem saying that $\Grn(\C)=\bigsqcup_{n\geq 0}\;
CM_n,$ is a
 union of the Calogero-Moser spaces. Note that our proof is purely algebraic
(as opposed to \cite{W}) and totally different from that in
\cite{W}.
\medskip

This paper is organized as follows. In \S3 we prove Theorem
\ref{DRiso} by means of  a $\B$-module version of
Kashawara's Theorem describing ${\mathcal{D}}$-modules 
concentrated on a point. In \S4 we re-interpret
$\B$-modules in terms of noncommutative geometry,
 and prove Theorem
\ref{jstariso} and  Theorem \ref{symb-restr}. Sections
5 and 6 contain proofs of Theorems \ref{framedsheaves} and
\ref{prb},
respectively. Appendix~A deals with modifications that one has to introduce
in the formalism of \cite{AZ} in order to be able to work with polygraded
algebras.
Finally, in Appendix B we prove some technical results on the
``noncommutative
surface" $\PPP,$ including Serre Duality and Beilinson Spectral Sequence.

\section{Kashiwara Theorem and De Rham Functor}
In this section we prove Theorem \ref{DRiso} by reducing it to a 
deformed version of Kashiwara's theorem on $\mathcal{D}$-modules
supported on a single point. 

To begin the proof of  Theorem \ref{DRiso} observe first
that any primary decomposable subspace 
$p\cd(W\gotimes\CG[x])\subset \CU \subset\frac1p\cd(W\gotimes\CG[x])$
is determined by the subspace\linebreak
$\frac{p\cdot \CU}{p^2\cdot(W\gotimes\CG[x])}$
$\subset
\frac{W\gotimes\CG[x]}{p^2\cdot(W\gotimes\CG[x])}$.
Similarly, any fat $\B$-submodule
$p\cdot(W\gotimes\B)\subset N \subset\frac1p\cd(W\gotimes\B)$
is determined by the 
 $\B$-submodule 
$\frac{p\cdot N}{p^2\cdot(W\gotimes\B)}\subset \frac{W\gotimes\B}{p^2\cdot
(W\gotimes\B)}$.
Observe further that the De Rham functor $\DR$ is right exact and
the homological dimension of the category of $\B$-modules equals~1
(see~\cite{CBH}).
Moreover, since $W\gotimes\B$ and $p^2\cd(W\gotimes\B)$
are projective $\B$-modules we get
$$
\DR\left(\frac{p\cd N}{p^2\cd(W\gotimes\B)}\right) = 
\frac{p\cd\DR(N)}{p^2\cd(W\gotimes\CG[x])}\; \subset\;
\frac{W\gotimes\CG[x]}{p^2\cd(W\gotimes\CG[x])}\;.
% = 
%\DR((W\gotimes\B)/p^2(W\gotimes\B)).
$$
Therefore, to prove Theorem~\ref{DRiso} it suffices, 
according to the definitions of $\Gradb (W)$ and $\Gr(W)$, to
show that the functor $\DR$ induces a bijection
between the following sets

{(a) the set of $\B$-submodules 
of $(W \gotimes \B)/ p^2\cdot (W \gotimes  \B)$;
and

(b) the set of vector subspaces in 
$(W \gotimes \C \Gamma [x])/p^2\cd (W \gotimes\C\Gamma[x])$ 
which are compatible with the direct sum decomposition}
$$
\frac{W\gotimes\CG[x]}{ p^2\cdot(W\gotimes \CG[x])}\;\; =\;\;
\bigoplus_{\mu}\;\;
\frac{W\gotimes\CG[x]}{(x-\mu)^{2 k_\mu}W\gotimes\CG[x]}
$$
where $p(x) = \prod_{\mu}(x-\mu)^{k_\mu}$ is a fixed $\Gamma$-invariant
polynomial, and where
we have used an identification: $
\frac{1}{p}\cdot(W \gotimes  \B)/ p\cdot (W \gotimes \B)\iso
(W \gotimes \B)/p^2\cdot (W \gotimes \B)$
provided by multiplication by $p$ (and similarly for $\CG[x]$).
\bigskip

Next, equip the vector 
space $(W\gotimes \CG[x])/p^2\cd (W\gotimes\CG[x])$ with a new
$\C\Gamma[x]$-module structure by requiring that $x \in \C\Gamma[x]$ acts on
$(W\gotimes\CG[x])\big/\bigl((x-\mu)^{2k_\mu}\cd W\gotimes\CG[x]\bigr)$
 as multiplication by
$\mu$, and $\Gamma$
acts as before. In other words, we replace the natural $x$-action by its
semisimple part. Let $\SS(W, p^2)$
denote the  result of such a
 semisimplification. 
The key observation is that a vector subspace
$\CU\subset (W\gotimes\CG[x])/p^2(W\gotimes\CG[x])$
is compatible with the direct sum decomposition as in (b) above,
if and only if $\CU$ is a $\CG[x]$-submodule in $\SS(W,p^2)$.
Thus Theorem~\ref{DRiso} reduces to the assertion 
that $\DR$ induces a bijection between

{(a) the set of $\B$-submodules of $(W \gotimes \B)/ p^2\cdot (W \gotimes \B)$;
and

(b) the set of $\C\Gamma[x]$-submodules in $\SS(W, p^2)$.}
\bigskip   

Our next step is to show that
the polynomial $p^2=p(x)^2$ above can be replaced by a simpler
polynomial.
For any $\mu \in \C,$ let $\Gamma_{\mu}$ be the stabilizer of $\mu$ in
$\Gamma$,
$m_{\mu}$ the order of $\Gamma/\Gamma_{\mu}$ and~$p_{\mu}$ the minimal
$\Gamma$-semiinvariant
polynomial vanishing on $\Gamma \cdot \mu$. In other words, we put
$$
\begin{array}{lllll}
\Gamma_{\mu} := \{1\}, & m_{\mu} := m,  & p_{\mu}(x) := x^m - \mu^m &&
\textrm{if  } \mu \neq 0;
\textrm{  and} \\ \Gamma_{\mu} := \Gamma, & m_{\mu} := 1, & p_{\mu}(x) := x &&
\textrm{if  } \mu = 0.
\end{array}
$$
Then any $\G$-invariant polynomial $p(x)$ can be factored as
$\dis\,
p(x) = \prod\nolimits_{\mu \in \C/\G}\; p_{\mu}(x)^{s_{\mu}}
\,.\,$
This factorization induces direct sum decompositions:
$$
\frac{W \gotimes \B}{p^2\cd (W \gotimes  \B)} =
\bigoplus_{\mu \in \C/\G} \;
\frac{W \gotimes \B}{p_{\mu}^{\;2 s_{\mu}}\cd(W \gotimes \B)}; \qquad
\SS(W, p^2) = \bigoplus_{\mu \in \C/\G}\; \SS(W, p_{\mu}^{2s_{\mu}})\;\;.
$$

The following is clear
\begin{lemma}
For any $\B$-submodule $N\subset \bigoplus_{\mu \in \C/\G}\;
\frac{W \gotimes \B}{p_{\mu}^{2 s_{\mu}}\cdot(W \gotimes \B)},$
we have
$$\qquad N\;\; =\;\; \bigoplus\nolimits_{\mu \in \C/\G}\;\;
\left(N\;\bigcap\;
\frac{W \gotimes \B}{p_{\mu}^{2 s_{\mu}}\cd(W \gotimes \B)}\right)\;.
{\hskip 30mm}\square
$$
\end{lemma}

Due to the above Lemma we may (and will)
assume without any loss of generality 
that $p(x) = p_{\mu}(x)^{2s_{\mu}},$ for some fixed
$\mu \in \C$ and some $s_{\mu}=1,2,\ldots$.
Further, we have
 $\SS(W, p_{\mu}^{\;2 s_{\mu}}) \simeq \SS(W, p_{\mu})^{\oplus
2 s_{\mu}}$, and this space is,  in effect,
a module over the quotient algebra $\C\G[x]/p_{\mu} \C\G[x]$. The
set of submodules in
$\SS(W, p_{\mu}^{\;2 s_{\mu}})$ may be therefore described by
the following result, which is proved by a straightforward computation.

\begin{lemma} \label{morita} (a) The correspondence $\CU \mapsto \CU
\otimes_{\Gamma_{\mu}} \Big(\C\G[x]/(x-\mu)
\C\G[x] \Big)$ establishes a Morita equivalence between the category
$\Rep(\Gamma_{\mu})$ of finite-dimensional representations of $\Gamma_{\mu}$
and the category of
finite-dimensional $\C\G[x]\big/p_{\mu} \C\G[x]$-modules.

(b) The $\Gamma_{\mu}$-module $\CU(W, p_{\mu}^{\;2 s_{\mu}})$ corresponding
to
$\SS(W, p_{\mu}^{\;2 s_{\mu}})$ via this equivalence, is equal to 
$W^{\oplus 2 s_{\mu}},$
viewed as a vector space (=module over $\Gamma_{\mu} = \{1\}$)
if $\mu \neq 0$ and as a
$\Gamma$-module $W \otimes_{\C} \bigl(\C[x]/x^{2 s_0}
\C[x]\bigr),$
if $\mu = 0$. \hfill $\square$
\end{lemma}

Thus, to prove Theorem \ref{DRiso} we have to establish a bijection
between 
the following sets

(a) The set of all $\B$-submodules of $W \gotimes \B / 
p_{\mu}^{\;2s_{\mu}}\cd(W \gotimes\B)$;\quad and

(b) The set of all $\Gamma_{\mu}$-submodules of $\CU(W, p_{\mu}^{\,2
s_{\mu}})$.
\medskip

\noindent
To that end, we introduce the following
\begin{definition}
Denote by $\mmod_\Gmu(\B)$ the category of all finitely generated
$\B$-modules $\CM$, such that
$p_\mu(x)$ acts locally nilpotently on $\CM$.
\end{definition}

If $\CM$ is an object in $\mmod_\Gmu(\B)$ then
$$
K_\mu(\CM)\,: =\, \Ker (x - \mu) \subset \CM
$$
is a $\Gamma_\mu$-module. Moreover, it is clear that the assignment
$\CM\mapsto K_\mu(\CM)$ gives  a functor
$K_\mu: \mmod_\Gmu(\B) \to \Rep(\Gamma_\mu)$.  Further, consider the induction
functor
$$
I_\mu:\Rep(\Gamma_\mu)\to\mmod_\Gmu(\B),\qquad
\CU \mapsto \CU\otimes_{\CG_\mu[x]}\B,
$$
where $\CG_\mu[x] = \CG_\mu\otimes\C[x]$ and 
the $\CG_\mu[x]$-module structure on $\CU\in\Rep(\Gamma_\mu)$
is given by the standard action of $\Gamma_\mu$ and the action of $x$
by the $\mu$-multiplication.

The following theorem is a deformed (and equivariant)
analogue of a well-known result of Kashiwara saying that any
$\D$-module concentrated at a point is the $\D$-module direct image of
a vector space (= $\D$-module on that point).

\begin{theorem}[Kashiwara Theorem]\label{kasha} Assume that the 
element $\tau \in \C\G$
involved in the definition of $\B$ is generic in the sense of Definition
\ref{def_gen}. Then
the functors $K_\mu$ and $I_\mu$ give mutually inverse equivalences 
between the
categories $\mmod_\Gmu(\B)$ and $\Rep(\Gamma_\mu)$.
\end{theorem}

Before we prove this theorem we record
a few  consequences of the condition: $\tau$
{\it is generic}. For any $k=1,2,\ldots,$ and any
integers $0\leq a\leq b$,
we define elements $\tau^{(k)}\,,\,\tau_{[a,b]}\in \CG$ by the equations
$$
y^k\cdot\tau=\tau^{(k)}\cdot y^k,\qquad
\tau_{[a,b]} = \sum\nolimits_{k=a}^b\;\tau^{(k)}\quad.
$$
The definition yields
\begin{equation}\label{tcomm}
y\cdot \tau_{[a,b]}=\tau_{[a+1,b+1]}\cdot y,\qquad\text{and}\qquad
x\cdot \tau_{[a,b]}=\tau_{[a-1,b-1]}\cdot x.
\end{equation}

\begin{lemma}\label{gentau}
(a)
The element $\tau$ is generic
if and only if  for all $a\le b$ the element $\tau_{[a,b]}\in \CG$ is
invertible. Furthermore, in this case for any $a\in\Z$ the element
$\tau_{[a,a+m-1]}$ acts by a constant {\rm(}independent of $a${\rm)}
in any representation of $\Gamma$.

(b) If $\tau$ is generic then for any $\mu\in\C$ and all
$a\le b$ the element $\sum_{i=a}^b\tau_{[i,i+m_\mu-1]}\in\CG$
is invertible.

(c) We have $[y\,,\,p_\mu(x)] = \frac{\tau_{[0,m_\mu-1]}}{m_\mu}
\cd p'_\mu(x)$.
\end{lemma}
{\sl Proof:}  To prove (a) recall (cf. e.g. [CBH]) that 
McKay correspondence associates to the cyclic group $\Z/m\Z=\bmu$
the affine Dynkin graph ${\widetilde A}_{m-1}$. Using an explicit
expression for  the roots it is easy to deduce the assertion.

To prove (b) note that if $\mu\ne 0$ then $m_\mu=m$ and the sum
in
question equals $|\tau|\cdot(b-a+1),$ where $|\tau|$ is the constant of part (a).
Hence this sum is invertible. If $\mu=0$ then $m_\mu=1$ and 
we have: $\sum_{i=a}^b\tau_{[i,i+m_\mu-1]}=\tau_{[a,b]}$.
Hence this
element is  invertible also.

Part (c) is proved by a direct computation. \hfill\qed

\bigskip
\noindent
{\sl Proof of Theorem \ref{kasha}:} It follows from the definitions of
$K_{\mu}$ and
$I_{\mu}$ that the functor $K_{\mu}$ is right adjoint to $I_{\mu}$. This
gives canonical adjunction
morphisms $\CU \to K_{\mu}(I_{\mu}(\CU))$ and $I_{\mu}(K_{\mu}(\CM)) \to
\CM$ for
any $\CM \in Ob(\mmod_\Gmu \B)$ and $\CU \in Ob (\Rep \Gamma_{\mu})$.
\bigskip

To show that $I_{\mu}(K_{\mu}(\CM)) \to \CM$ is an isomorphism,
 set $\CM_k \,:=\,
\Ker p_\mu(x)^{k+1}
\subset \CM$ and write $M_k:=\CM_k/\CM_{k-1},$ for
short. On $\oplus_k\, M_k,$
we have the following structure.

First, it is clear that the increasing
filtration $\,\{\CM_k\}_{k=0,1,\ldots}\,$ is stable under the
action of the subalgebra $\Gamma[x]\subset\B$, hence each $M_k$ is a
$\Gamma[x]$-module. Further, multiplication by $p_\mu(x)$
takes $\CM_k$ to $\CM_{k-1}$ and thus induces a map
$p:M_k\to M_{k-1}$. Moreover, the action of $y$ moves $\CM_k$ to
$\CM_{k+1}$ and thus induces a map $M_k\to M_{k+1}$. Finally, it is
clear that the map $p$ is an embedding of $\Gamma[x]$-modules, while
$y$ is a morphism of $\C[x]$-modules, and $p'_\mu(x):M_k\to M_k$
is an isomorphism (because $p_\mu(x)$ and $p'_\mu(x)$ are coprime).
Let us prove by induction in $k$ that
\begin{equation}\label{pym} 
(p\cdot y)|_{M_k} =
- \sum_{i=0}^k\frac{\tau_{[i,i+m_\mu-1]}}{m_\mu}\cdot p'_\mu(x),\quad
(y\cdot p)|_{M_{k+1}} =
- \sum_{i=1}^{k+1}\frac{\tau_{[i,i+m_\mu-1]}}{m_\mu}\cdot p'_\mu(x).
\end{equation}
The base of induction ($k=-1$) is clear. Assume that we have
verified~(\ref{pym}) for $k-1$. Then for any $a\in M_k,$ applying
Lemma \ref{gentau}(c) and the induction hypothesis, we get
$$
p\cdot y\cdot a = y\cdot p\cdot a - [y,p]\cdot a =
- \sum_{i=1}^k\frac{\tau_{[i,i+m_\mu-1]}}{m_\mu}\cdot p'_\mu a
- \frac{\tau_{[0,m_\mu-1]}}{m_\mu}\cdot p'_\mu a =
- \sum_{i=0}^k\frac{\tau_{[i,i+m_\mu-1]}}{m_\mu}\cdot p'_\mu a.
$$
Note that since $\tau$ is generic it follows from Lemma \ref{gentau}(b) that
the map 
$p\cd y:M_k\to M_k$ is a bijection. On the other hand, $p$ is 
injective by definition. Hence $y$ gives an isomorphism $M_k\iso
M_{k+1}$. It follows that, for any $b\in M_{k+1}$
there exists an $a\in M_k$ such that  $b=y\cd a$.
Applying Lemma \ref{gentau}(b), we calculate
$$
\begin{array}{l}
y\cdot p\cdot b = y\cdot p\cdot y\cdot a =
-y\cdot  \sum_{i=0}^k\frac{\tau_{[i,i+m_\mu-1]}}{m_\mu}\cdot p'_\mu\cdot a =
-\sum_{i=0}^k\frac{\tau_{[i+1,i+m_\mu]}}{m_\mu}\cdot y\cdot p'_\mu\cdot a =
\medskip\\\qquad\qquad\qquad\qquad\qquad
-\sum_{i=1}^{k+1}\frac{\tau_{[i,i+m_\mu-1]}}{m_\mu}\cdot p'_\mu\cdot y\cdot a =
- \sum_{i=1}^{k+1}\frac{\tau_{[i,i+m_\mu-1]}}{m_\mu}\cdot p'_\mu\cdot b
\end{array}
$$
(In the third equality
we use the fact that the operators 
$y:M_k\to M_{k+1}$ and $p'_\mu(x):M_k\to M_k$ commute, because
their commutator on $\CM$ takes $\CM_k$ to $\CM_k$, hence
induces the zero map $M_k\to M_{k+1}$).
Thus we have proved~(\ref{pym}) for any $k$ and, moreover,
we have shown along the way that for any $k$ the map
$y:M_k\to M_{k+1}$ is an isomorphism. This means that the action of
the subalgebra $\C[y]\subset\B$ gives an isomorphism
$\,
\CM \cong \CM_0\otimes\C[y].
$ Further, we have:
 $\CM_0 = \Ker p_{\mu}(x)$ and $K_{\mu} (\CM) = \Ker (x - \mu)$.
We use 
the equalities:
$$
\Ker p_\mu(x) =
\!\bigoplus_{\gamma\in\Gamma/\Gamma_\mu}\!
\Ker(\gamma(x-\mu)\gamma^{-1}) =
\!\bigoplus_{\gamma\in\Gamma/\Gamma_\mu}\!
\left(\Ker(x-\mu)\right)\gamma^{-1} =
\Ker(x-\mu)\otimes_{\C[\Gamma_\mu]}\CG
$$
to conclude that $I_{\mu}(K_{\mu}(\CM)) \to \CM $ is an isomorphism.

\bigskip

To show that the canonical morphism $f: \CU\to K_\mu(I_\mu(\CU))$ is an
isomorphism,
first note that it is clearly injective. Hence we have an exact
sequence
$$
0 \too \CU \too K_\mu(I_\mu(\CU)) \too \CU' \too 0
\quad\text{where}\quad \CU':= {\mathbf{Coker}}(f)\,.
$$
On the other hand, the functor $I_\mu$ is exact since $\B$ is a
flat $\Gamma_\mu[x]$-module. Hence,
applying the functor $I_\mu(-)$, we obtain an exact sequence
$$
0 \too I_\mu(\CU) \stackrel\alpha\too 
I_\mu(K_\mu(I_\mu(\CU))) \too I_\mu(\CU') \too 0.
$$
The argument of the first part of the proof applied to the
$\B$-module
$I_\mu(\CU)$, shows
that the morphism $\alpha$ above is an isomorphism. Hence
$I_\mu(\CU')=0$. But this clearly yields $\CU'=0$. Thus,
the map $\CU\to K_\mu(I_\mu(\CU))$ is an isomorphism,
and  Theorem
\ref{kasha} follows. \hfill\qed

\bigskip
\noindent
{\sl End of proof of Theorem \ref{DRiso}:} 
The De Rham functor $\DR$ restricted to the set of submodules 
in $(W \gotimes \B)/p_{\mu}^{\;2 s_{\mu}}\cd  (W \gotimes \B)$ 
can be factored as a composition of the equivalence 
$K_{\mu}: \mmod_\Gmu \B \to \Rep \Gamma_{\mu},$
and the Morita equivalence 
$\Rep \Gamma_{\mu} \to {\mathop{{\sf mod}}} \;
\bigl(\C\G[x]/p_{\mu}\cd\C\G[x] \bigr)$
of Lemma \ref{morita}~(b). 
Hence it is an equivalence as well. By a straightforward 
(but a bit tedious) computation one deduces that
$$ 
K_{\mu}\bigl((W \gotimes \B)/p_{\mu}^{\;2 s_{\mu}}\cd  (W \gotimes
\B)\bigr)
 \simeq
\CU(W, p_{\mu}^{\;2 s_{\mu}})
$$
where $\CU(W, p_{\mu}^{\;2 s_{\mu}})$ is given by Lemma \ref{morita} (b).
This implies Theorem \ref{DRiso}, as we have seen
in the first half of this section.

Finally, one can check that the map $\Diff$ defined just below
the statement of Theorem~\ref{DRiso} is in effect  the inverse bijection:
$\Gr(W)\to\Gradb(W)$.
\hfill\qed

\section{$\B$-module Grassmanian and Sheaves on $\PPP$}
Recall the
bigraded algebra
$\Qu$ defined in \eqref{qu}.
Let $\gr^2(\Qu)$ be the
category of finitely generated bigraded right $\Qu$-modules $M=\oplus
M_{i,j}$.
Let
$\tor^2(\Qu)$ denote
its Serre subcategory formed by all modules $M$ such that there
exists a pair 
$(i_0,j_0)$ such that for any $i>i_0$ and $j>j_0$ we have $M_{i,j}=0$.
Consider the quotient category $\qgr^2(\Qu) = \gr^2(\Qu)/\tor^2(\Qu)$.
The category $\qgr^2(\Qu)$
will be viewed as the category of coherent sheaves on a noncommutative
scheme  $\PPP$
(see Appendix A for details). Thus,
\emph{by definition}, we put $\coh(\PPP):=\qgr^2(\Qu),$
and we write
$\pi:\gr^2(\Qu)\to\coh(\PPP)$
for the canonical projection functor.

The isomorphism $\B\cong \Qu/((z-1)\Qu + (w-1)\Qu)$ gives rise to
a `restriction'  functor
$$
j^*:\coh(\PPP)\to\mmod(\B),\qquad
E = \pi(M) \longmapsto \dlimkl M_{k,l}\cong M/((z-1)M + (w-1)M),
$$
where the direct limit is taken with respect to the maps
$M_{k,l}\to M_{k+1,l}$ and $M_{k,l}\to M_{k,l+1},$
induced by  multiplication by $z$ and $w$, respectively.

There are canonical isomorphisms
$$
\Qu/z\Qu\; \cong\; (\C[x]\otimes\C[y,w])\#\Gamma,\qquad
\Qu/w\Qu \;\cong\; (\C[x,z]\otimes\C[y])\#\Gamma.
$$
Thus, we obtain the following equivalences of categories: 
$\qgr^2(\Qu/z\Qu) \simeq \qgr(\C[xy,xw]\#\Gamma)$, and
$\qgr^2(\Qu/w\Qu)\,\simeq\,\qgr(\C[xy,yz]\#\Gamma)$,
(see Appendix, Corollary~\ref{qgrrqgr}).
The two categories on the right can be viewed as the
categories of  $\Gamma$-equivariant 
coherent sheaves on the ordinary
projective line $\PP^1$. We denote the corresponding copies
of  $\PP^1$  by $\PP^1_z$
and $\PP^1_w$ respectively. We have
the corresponding push forward and pull back functors:
$$
\begin{array}{rr}
\dis(\iz)_*:\coh(\PP^1_z)\to\coh(\PPP), &
\iz^*:\coh(\PPP)\to\coh(\PP^1_z),\\
\dis(\iw)_*:\coh(\PP^1_w)\to\coh(\PPP), &
\iw^*:\coh(\PPP)\to\coh(\PP^1_w),
\end{array}
$$
Let $L^1i^*$ stand for the first derive functor.

Given  a $\B$-module $\CM$,
 we let its {\em support} be the support of  $\CM$,
viewed as a module
over the subalgebra $\C[x]\subset \B$
(more precisely, the union of supports of
all $\C[x]$-finitely generated
submodules in $\CM$).

\begin{definition} \vi
Let $\Quot(W)$ be the category of all surjections $W\gotimes\CO\onto F$
in the category $\coh(\PPP)$ such that
\begin{equation}\label{iziwf}
\iz^*F=L^1\iz^*F=0\quad\text{and}\quad
L^1\iw^*F=0.
\end{equation}

\vii
Let $\BQuot(W)$ be the category of all  $\B$-module
surjections $W\gotimes\B\onto\CM$,
 such that $\CM$ has zero-dimensional support.
\end{definition}

\begin{theorem}\label{jstarth}
The functor $j^*$ takes any object of the category $\Quot(W)$
to an object of  $\BQuot(W)$ and, moreover, gives an
equivalence $\Quot(W)\iso\BQuot(W)$.
\end{theorem}
{\sl Proof:}
First note that $j^*$ is exact, being a  direct limit functor.
Thus, to show that $j^*$ takes $\Quot(W)$ to $\BQuot(W)$ it suffices to show
that for any object of $\Quot(W)$ of the form
$W\otimes\CO\onto F$, the $\B$-module
$j^*F$ has zero-dimensional support. Indeed, let $F=\pi(M)$ where $M = \oplus
M_{k,l}$.
Then conditions $\iz^*F=L^1\iz^*F=0$ imply that 
$z$-multiplication gives an isomorphism $M_{k,l}\iso M_{k+1,l}$
for $k, l \gg 0$. Since $M$ is finitely
generated we can choose $k_0$ and $l_0$ such that
$\oplus_{k\ge k_0,\ l\ge l_0}M_{k,l}$ is generated by $M_{k_0,l_0}$
and   such that
$z: M_{k,l}\to M_{k+1,l}$ is an isomorphism for all $k\ge k_0$, $l\ge
l_0$.
Let $p$ be the characteristic polynomial of the operator
$z^{-1}x:M_{k_0,l_0} \to M_{k_0,l_0}$. Then it is easy to see
that $p(x)$ acts locally nilpotently on $j^*F = \dlim M_{k,l}$
and hence $j^*$ indeed defines a functor $j^*: \Quot(W) \to \BQuot(W)$.

The assertion that $j^*$ is an equivalence will be proved by constructing a
quasi-inverse
functor $j_{!*}$. To that end, let $\psi:W\gotimes\B\onto\CM$ be an
object of $\BQuot(W)$. Since $\CM$ has zero-dimensional support
there exists a $\Gamma$-invariant polynomial $p(x)$ that acts by
zero on the subspace $\psi(W\gotimes\CG)\subset\CM$. It is clear
that we have $\psi(W\otimes p(x)) = 0$, hence
\begin{equation}\label{annihilates}
\psi(p\cd (W\gotimes \B)) = 0
\end{equation}
Let $\B_{k,l}$ be the natural increasing bifiltration of $\B$ (induced by the
bigrading of $\Qu$) and
$$
\CM_{k,l} = \psi(W\gotimes\B_{k,l}) \subset \CM,
$$
 the induced bifiltration of $\CM$. It follows from (\ref{annihilates})
that this bifiltration stabilizes with respect to the first index
when $k\ge d = \deg p(x)$, that is we have
$$
\CM_{k,l} = \psi(W\gotimes\B_{k,l}) =
\psi(W\gotimes\B_{d,l}) = \CM_{d,l} \subset \CM
\qquad\text{for $k\ge d$ and all $l$.}
$$
It follows from the definition that the bifiltration $\CM_{k,l}$ is
compatible with the bifiltration on~$\B$. Moreover, it is clearly
increasing, finitely generated and exhaustive (because $\psi$
is surjective). Hence $M = \oplus_{k,l}\CM_{k,l}$ is a finitely
generated $\Qu$-module, where the action  on $M$ of $x$-generators and $y$-generators
of $\Qu$ is given by the $x$ and $y$ multiplication maps 
$\CM_{k,l}\to\CM_{k+1,l}$ and $\CM_{k,l}\to\CM_{k,l+1}$ respectively,
and the action of $z$ and $w$ generators is given by the natural embeddings
$\CM_{k,l}\into\CM_{k+1,l}$ and $\CM_{k,l}\into\CM_{k,l+1}$ respectively.

Consider  $F = \pi(M)$, a coherent sheaf  on $\PPP$. It follows from
the definition of~$F$ that the $z$-multiplication map $F\to F(1,0)$
is an isomorphism (since $\CM_{k,l}=\CM_{k+1,l}$ for $k\ge d$), hence
$\iz^*F=L^1\iz^*F=0$. On the other hand, the $w$-multiplication map
$F\to F(0,1)$ is an embedding (because $\CM_{k,l}\subset\CM_{k,l+1}$),
hence $L^1\iw^*F=0$. Finally note that the map $\psi$ is compatible
with the bifiltrations on $W\gotimes\B$ and $\CM$, hence it gives rise
to a map $\tilde\psi:W\gotimes\CO\to F$ of coherent sheaves on $\PPP$.
Moreover, this map is surjective by definition, hence it gives an object of
$\Quot(W)$. Finally, it is easy to show that this way
we obtain a functor $j_{!*}:\BQuot(W)\to\Quot(W)$.

We show that $j_{!*}$ and $j^*$ are quasi-inverse. Let
$W\gotimes\CO\stackrel\psi\to F$ be an object of the category $\Quot(W)$
and let $M = \oplus M_{k,l}$, where $M_{k,l}=H^0(\PPP,F(k,l))$.
Then $M$ is a bigraded $\Qu$-module and it is clear that $\pi(M) = F$.
Note that we have exact sequences
$$
0 \to {\iz}_*L^1\iz^*F(k+1,l) \to F(k,l) \to F(k+1,l) 
\to {\iz}_*\iz^*F(k+1,l) \to 0,
$$
$$
0 \to {\iw}_*L^1\iw^*F(k,l+1) \to F(k,l) \to F(k,l+1) 
\to {\iw}_*\iw^*F(k,l+1) \to 0.
$$
Moreover, applying~(\ref{iziwf}) we get
$$
\begin{array}{l}
L^1\iz^*F(k+1,l) \cong (L^1\iz^*F)(l) = 0,\qquad
\iz^*F(k+1,l) \cong (\iz^*F)(l) = 0,\\
L^1\iw^*F(k,l+1) \cong (L^1\iw^*F)(k) = 0.
\end{array}
$$
Combining these isomorphisms with the above exact sequences and with the
definition
of $M_{k,l}$, we see that the maps $M_{k,l}\stackrel{z}\to M_{k+1,l}$ and
$M_{k,l}\stackrel{w}\to M_{k,l+1}$ are an isomorphism and an embedding
respectively. Therefore $j^*F=\bigcup_l M_{k,l}$ for any $k$. On the
other hand we have $\B_{k,l}=H^0(\PPP,\CO(k,l))$ by
definition and it is clear that the map $j^*\psi$ sends
$W\gotimes\B_{k,l}$ to $M_{k,l}$ and coincides there with the map
$H^0(\psi(k,l))$. Thus to show that $j_{!*}(j^*F)\cong F$ it suffices to
show that this map is surjective for all $k$ and $l$ sufficiently
large. But this is nothing but the definition of the map $\psi$ being a
surjection in the category $\coh(\PPP)$.

Now assume that $W\gotimes\B\to\CM$ is an object of $\BQuot(W)$. Then
by definition of $j_{!*}$ we have $\CM_{k+1,l}=\CM_{k,l}$ for $k\ge d$
and all $l$, and when $k\ge d$ is fixed the filtration $\CM_{k,l}$ of
$\CM$ is exhaustive. Hence $\CM=\dlim\CM_{k,l}$, that is
$j^*j_{!*}(\CM)=\CM$.
\hfill\qed

\bigskip

Now we give a more rigorous version of Definition \ref{nc-pre}.
Let $E$ be a coherent sheaf on $\PPP$ such that
$\iz^*E\cong W\gotimes\CO_{\PP^1_z}$. Recall that
$\eps$ denotes a fixed primitive character: $\G=\bmu\into \C^*$.

\begin{definition}\label{trivdef1} \vi
We say that the sheaf $E$ is \emph{trivialized in a neighborhood
of $\PP^1_z$} if we are given embeddings
$$
(W\otimes \eps^{-n})\gotimes\CO(-n,0) \stackrel\phi\into E
\stackrel\psi\into (W\otimes \eps^{n})\gotimes\CO(n,0)
$$
such that the composite
$$
\psi\ccirc\phi\in
\Hom((W\otimes \eps^{-n})\gotimes\CO(-n,0),(W\otimes \eps^{n})\gotimes\CO(n,0)) \cong
\Hom_\Gamma\bigl(W,(W\otimes \eps^{2n})\gotimes \Qu_{2n,0}\bigr)
$$
equals multiplication by $P(x,z)^2$, where  
$P(x,z)\in\C[x,z]$ is a $\Gamma$-semiinvariant homogeneous 
polynomial of degree $n$, such that $P(1,0)=1$.

\vii
We call two trivializations $(\phi,\psi)$ and $(\phi',\psi')$
of the sheaf $E$  \emph{equivalent} if there exists a pair of
$\Gamma$-semiinvariant homogeneous polynomials $\tq(x,z)$ and
$\tq'(x,z)$ such that $\tq(1,0) = \tq'(1,0) = 1$ 
and the following diagram commutes
$$
\xymatrix{
(W\otimes \eps^{-n})\gotimes\CO(-n,0) \ar[r]_-\phi &
E \ar[r]_-\psi \ar@{=}[dd] &
(W\otimes \eps^{n})\gotimes\CO(n,0) \ar[d]_{\tq}
\\
(W\otimes \eps^{-n''})\gotimes\CO(-n'',0) \ar[u]_{\tq} \ar[d]^{\tq'}
&& 
(W\otimes \eps^{n''})\gotimes\CO(n'',0) 
\\
(W\otimes \eps^{-n'})\gotimes\CO(-n',0) \ar[r]^-{\phi'} &
E \ar[r]^-{\psi'} &
(W\otimes \eps^{n'})\gotimes\CO(n',0) \ar[u]^{\tq'}\;\;.
}
$$
\end{definition}

\begin{remark}
We can always replace a trivialization by an equivalent one
with $n\equiv 0\, {\sl{mod}}\; m$, thus getting rid of $\eps^n$ and $\eps^{-n}$
factors in the definition and making the polynomial $P(x,z)$ 
$\Gamma$-invariant.
\end{remark}

\begin{definition}\label{trivdef2} Let $\Gradp(W)$
 be the set of all equivalence
classes of trivializations in a neighborhood of $\PP^1_z$
of coherent sheaves $E$ on $\PPP$ such that
$\iz^* E\cong W\gotimes\CO_{\PP^1_z}$.
\end{definition}

\bigskip
\noindent
{\sl Proof of Theorem \ref{jstariso}} ({\it bijection between $\Gradb(W)$ and
$\Gradp(W)$}):\;
For any $\Gamma$-invariant polynomial $p(x)=\sum_{k=0}^d\,a_kx^k,$
let $\Gradpp(W)\subset\Gradp(W)$
denote the set of all sheaves admitting a trivialization $(\phi,\psi)$
with $\psi\ccirc\phi=P(x,z)^2$, where $P(x,z)$ is the homogenization
of $p(x)$, that is $P(x,z)=\sum_{k=0}^d\,a_kx^kz^{d-k}.$ Then we have
$$
\Gradp(W) = \bigcup\nolimits_{p(x)}\;\Gradpp(W).
$$
We will show that the functor $j^*$ induces a bijection
between $\Gradpp(W) \subset \Gradp(W)$ and the subset 
$\Gradbp(W) \subset \Gradb(W)$ formed by all $\B$-submodules 
(or, equivalently, quotient modules) of
$\po(W \gotimes \B)/p\cd (W \gotimes \B)$.

Let $\Quot_p(W) \subset \Quot(W)$ be  the subset formed
by
surjections: $W \gotimes \CO \onto F,$ which send the image of
the map: $W\gotimes\CO(-2n,0) \stackrel{P(x,z)^2}\lra W\gotimes\CO$
 to zero in $F$. We may identify the set  $\Gradpp(W)$ with 
 $\Quot_p(W)$ via the assignment
$$
{\scriptstyle
\Big\{W\gotimes\CO(-n,0) \stackrel\phi\too E
\stackrel\psi\too W\gotimes\CO(n,0)\Big\}}\quad\longmapsto
\quad
{\scriptstyle
\Big\{W\gotimes\CO \too
\Coker\Big(E(-n,0) \stackrel{\psi(-n,0)}\too
 W\gotimes\CO\Big)\Big\}
}\;.
$$

Hence, Theorem \ref{jstarth} implies that the functor $j^*$
provides an identification
of the set\linebreak
$\Quot_p(W)$
with the subset $\BQuot_p(W)\subset \BQuot(W)$ formed
by all surjections
$W \gotimes \B$
$ \onto \CM$ which send $p^2\cdot(W \gotimes \B)$ to zero (in $\CM$).
On the other hand, to any object 
$$p\cdot(W\gotimes \B) \stackrel\phi\lra N
\stackrel\psi\lra \mbox{$\frac{1}{p}$}\cd(W\gotimes \B)
$$ in $\Gradbp(W)$ we associate
the quotient $$W\gotimes\B \onto \Coker\Big(p\cd N \stackrel{p\psi}\lra
W\gotimes\B\Big).$$ 
This yields an identification of $\BQuot_p(W)$ with $\Gradbp(W)$.
Therefore,
we get a bijection $\Gradpp(W) \iso$ $ \Gradbp(W)$.

Note that,  for any polynomial
$p(x)$ dividing $q(x)$, the map
$j^*$ commutes with the natural embeddings
$\Gradpp(W)\into\Gradpq(W)$
and $\Gradbp(W)\into\Gradbq(W)$.
The assertion of Theorem~$\ref{jstariso}$ follows.\hfill\qed
\medskip

Finally we prove Theorem~\ref{symb-restr}. Recall that
the pull back functor $\iw^*$ takes any sheaf
trivialized in a neighborhood of $\PP^1_z$, to a sheaf on
$\PP^1_w$ trivialized in a neighborhood of the point $P$,
which is the same as a $\Gamma$-equivariant sheaf on
the ordinary projective line $\PP^1$ trivialized in a Zariski
neighborhood of the infinity. Thus $\iw^*$ induces a map
$\Gradp(W) \to \Graff(W)$.

\bigskip
\noindent
{\sl Proof of Theorem \ref{symb-restr}:}
The claim easily follows from the definitions of the  maps involved.
In more details,
 let $N\subset p\cd (W\gotimes \B)$ be a fat $\B$-submodule
and $M = p\cd (W\gotimes \B)$. Then it follows that
$j_{!*}(N) = \pi(\oplus N_{k,l})$, where
$N_{k,l} = N\cap(W\gotimes\B_{k,l})$ and $j_{!*}$ is the quasi-inverse
to $j^*$ defined in the proof of Theorem \ref{jstarth}.
Hence $\iw^*(j_{!*}(N)) = \pi(\oplus N_{k,l}/N_{k,l-1})$ and
the restriction of this sheaf to $\PP^1_w \smallsetminus \{\infty\}$
is isomorphic to $\dlimk\, N_{k,l}/N_{k,l-1} = \Symb(N)$.
This shows that the diagram~(\ref{pic}) is commutative.
Further, the maps in the upper row of the diagram are bijections
by Theorems~\ref{DRiso} and \ref{jstariso}. Further, the maps
$\DR$ and $j^*$ are evidently $\Gw$-equivariant, hence $\Diff$
and $j_{!*}$ are $\Gw$-equivariant as well. And finally, the map
$\Symb$ is evidently $GL_{\CG(x)}(W\gotimes\CG(x))$-equivariant,
hence the maps $i_w^*$ and $\sigma$ are
$GL_{\CG(x)}(W\gotimes\CG(x))$-equivariant as well.
\qed

\section{Monads and Quiver Varieties}

The two lines $\PP^1_z$ and $\PP^1_w$ in $\PPP$ intersect at the
point $P$, corresponding to the quotient algebra
$$
\Qu/(z\Qu + w\Qu) \cong (\C[x]\otimes\C[y])\#\Gamma\;\;.
$$
Moreover, the category $\qgr^2(\Qu/(z\Qu+w\Qu))$ is equivalent
to the category of finite dimensional $\Gamma$-modules.
The point $P$ is given on the line $\PP^1_z$ by the equation
$w=0$ and on the line $\PP^1_w$ by the equation $z=0$.
Let $\izp:\{P\}\into\PP^1_w$ and $\iwp:\{P\}\into\PP^1_z$ 
denote the embeddings.
There is a canonical isomorphism of functors
$$
(\izp)^*\ccirc\iw^* \cong (\iwp)^*\ccirc\iz^* :\;
\coh(\PPP) \too \coh(P) = \Rep(\Gamma).
$$

Let $W, V$ be a pair of $\Gamma$-modules as in Definition \ref{quiverdef}.

\begin{definition} \label {quiverdef2} A  coherent sheaf $E$ on $\PPP$ 
is called 
$W$-framed  provided it is equipped
with  two isomorphisms $\iz^*E\cong W\gotimes\CO_{\PP^1_z}$ and
$\iw^*E\cong W\gotimes\CO_{\PP^1_w}$,
which agree at the point $P$.
\end{definition}

Let $\MPPP(V,W)$ denote the set of isomorphism classes of
$W$-framed torsion free sheaves $E$
(for the definition of `{\it torsion free}' see \cite[Def. 1.1.4]{BGK})
 on $\PPP$ such that
$H^1(\PPP,E(-1,-1))\cong V$.

\begin{theorem}\label{qs1}
The set $\MPPP(V,W)$ is in a natural bijection with
the quiver variety $\FM_\Gamma^\tau(V,W)$.
\end{theorem}
{\sl Sketch of  proof:}
The proof is essentially the same as that of  Theorem~1.3.10
of~\cite{BGK}, \S4. So, we will skip  most
of the details and only mention the points that are different from
\cite{BGK}, \S4.

The first difference is that, in the present situation, the monad
representing a framed sheaf  has a form
 slightly different from the one used in \cite{BGK}.
Specifically, for any point $(B_1,B_2,I,J)$
of the quiver variety, our monad is now
given by  the following complex
\begin{equation}\label{monadppp}
0 \to V\gotimes\CO(-1,-1)\; \stackrel a\lra\,
\begin{array}{c}
(V\otimes \eps)\gotimes\CO(0,-1) \\
\bigoplus \\
(V\otimes \eps^{-1})\gotimes\CO(-1,0) \\
\bigoplus \\
W\gotimes\CO
\end{array}\,
\stackrel b\lra\; V\gotimes \CO \to 0,
\end{equation}
\begin{equation}\label{monadppp1}
a = (  B_1z-x , B_2w-y, Jzw), \qquad
b = (-(B_2w-y), B_1z-x, I).
\end{equation}

Second, whenever the functor $i^*$ (the restriction
to the line at infinity in $\PP^2_\tau$) was used
in [BGK], it should be replaced by a pair
of functors $\iz^*$, $\iw^*$.

Third, (\cite{BGK}, Lemma~4.2.12) should be
 replaced by the following isomorphisms
\begin{align}\label{he}
& H^0(\PPP,E(-1,0)) = H^0(\PPP,E(0,-1)) = H^0(\PPP,E(-1,-1)) 
= 0,\nonumber\\
& H^2(\PPP,E) = H^2(\PPP,E(-1,0)) = \\
& \hp\hskip 30mm = 
H^2(\PPP,E(0,-1)) = H^2(\PPP,E(-1,-1)) = 0,\nonumber\\
& H^0(\PPP,E(-1,0)) = H^0(\PPP,E(0,-1)) = H^0(\PPP,E(-1,-1)),\nonumber
\end{align}
and furthermore there is a canonical exact sequence
\begin{equation}\label{he1}
0 \to H^0(\PPP,E) \to W \stackrel{f_E}\lra V \to H^1(\PPP,E) \to 0.
\end{equation}

Fourth, the Beilinson spectral sequence takes the
following form (see Appendix B)
\begin{equation}\label{BSSPPP}
E_1^{p,q} = 
\Big\{
{\scriptstyle
\Ext^q\bigl(\CO(1,1),E\bigr)\mathop{\otimes}\limits_\Gamma\CO(-1,-1) \to
\begin{array}{c}
\scriptstyle
\Bigl(\Ext^q\bigl(\CO(1,0),E\bigr)\otimes \eps\Bigr)
\mathop{\otimes}\limits_\Gamma\CO(-1,0)\\
\bigoplus \\
\scriptstyle
\Bigl(\Ext^q\bigl(\CO(0,1),E\bigr)\otimes \eps^{-1}\Bigr)
\mathop{\otimes}\limits_\Gamma\CO(0,-1)
\end{array}
\to \Ext^q(\CO,E)\mathop{\otimes}\limits_\Gamma\CO 
}
\Big\}
\end{equation}
%%%%%%%%%%%%%%%%%%%%%%%%%%%%%%%%%%%%%%%%%%%%%%%
%%\begin{equation}\label{BSSPPP}
%\begin{align}\label{BSSPPP}
%E_1^{p,q} = &\Big\{\Ext^q(\CO(1,1),E)\otimes\limits_\Gamma\CO(-1,-1)\; \too
%\,
%\begin{array}{c}
%(\Ext^q(\CO(1,0),E)\otimes \eps)\otimes\limits_\Gamma\CO(-1,0)\\
%\bigoplus \\
%(\Ext^q(\CO(0,1),E)\otimes \eps^{-1})\otimes\limits_\Gamma\CO(0,-1)
%\end{array}\nonumber\\
%&\too\;
%\Ext^q(\CO,E)\otimes\limits_\Gamma\CO \Big\}
%\end{align}
%%\end{equation}
%%%%%%%%%%%%%%%%%%%%%%%%%%%%%%%%%%%%%%%%%%%%%%%
We apply this spectral sequence to obtain a monadic description
of an arbitrary framed coherent sheaf $E$. Using~(\ref{he}) we see that,
for any $W$-framed sheaf $E$ such that\linebreak
$H^1(\PPP,E(-1,-1))\cong V$,
 the
spectral sequence  takes the form
$$
\xymatrix{
V\gotimes\CO(-1,-1)\; \ar[r] \ar@{-->}[drr]_{d_2^{-2,1}} &
\;(V\otimes \eps)\gotimes\CO(-1,0)\, \bigoplus\,
(V\otimes \eps^{-1})\gotimes\CO(0,-1)\; \ar[r] &
%(\Coker f)\gotimes\CO\\
%&& (\Ker f)\gotimes\CO
\;H^1(E)\gotimes\CO\\
&& H^0(E)\gotimes\CO
}
$$
Further, one  shows that one can replace $H^1(E)=\Coker f_E$ 
and $H^0(E)=\Ker f_E$ by $V$ and $W$ respectively and lift the 
differential $d_2^{-2,1}:E_2^{-2,1} \to E_2^{0,0}$ to 
a morphism
$V\gotimes\CO(-1,-1) \to W\gotimes\CO$. Finally,
replacing the spectral sequence with the total complex
one obtains the desired monadic description~(\ref{monadppp})
of the sheaf $E$. We leave for the reader  to check that
the maps in (\ref{monadppp})  take form (\ref{monadppp1})
for an appropriately chosen quiver data $(B_1,B_2,I,J)$.
\hfill\qed\medskip

\begin{remark}
There is an alternative way to prove Theorem~\ref{qs1}.
using the following trigraded algebra
$$
S := \C\langle\xi,\eta,\zeta,x,z,y,w\rangle\#\Gamma 
\left/\left\langle
\begin{array}{l}
{}[\bullet,\zeta] = [\bullet,z] = [\bullet,w] = [\xi,x] = [\eta,y] = 0,\\
{}[\eta,\xi]=\tau\zeta^2,\ [y,x]=\tau zw,\ 
{}[\eta,x] = \tau\zeta z,\ [y,\xi] = \tau\zeta w,\\
{}\xi z = x\zeta,\ \eta w = y\zeta
\end{array}
\right\rangle\right.
$$
Let $X$ be the corresponding noncommutative variety
(i.e., such that $\coh(X) = \qgr^3(S)$). Then we have
a diagram
\begin{equation}\label{radon}
\PP^2_\tau\;\stackrel{p}\longleftarrow\;X\;\stackrel{q}\longrightarrow\;
\PPP\;\;,
\end{equation}
where the morphism 
$q$ is a noncommutative analog of the blowup of the point 
$P$ on $\PPP$, and  the morphism $p$ is a noncommutative analog of
the blowup of a pair of points on the line at infinity.
One can show that  Fourier--Mukai type functors:
$$
q_*p^*:\coh(\PP^2_\tau) \to \coh(\PPP),\quad\text{and}\quad
p_*q^*:\coh(\PPP) \to \coh(\PP^2_\tau)
$$
induce mutually inverse bijections between the corresponding
sets of (isomorphism classes
of) $W$-framed torsion free sheaves.
Theorem~\ref{qs1} is now immediate from (\cite{BGK}, Theorem~1.3.10).
\end{remark}\medskip

We now turn to the proof of Theorem \ref{framedsheaves}.
Let $E$ be a $W$-framed torsion free coherent sheaf
on $\PPP$ such that  $H^1(\PPP,E(-1,-1))$ $\cong V$. Theorem~\ref{qs1}
implies that $E$ can be represented as the cohomology sheaf
of  monad~(\ref{monadppp}).
We consider the following maps
$$
(W\otimes \eps^{-n})\gotimes\CO(-n,0)\; \stackrel\Phi\lra\,
\begin{array}{c}
(V\otimes \eps)\gotimes\CO(0,-1) \\
\bigoplus \\
(V\otimes \eps^{-1})\gotimes\CO(-1,0) \\
\bigoplus \\
W\gotimes\CO
\end{array}\,
\stackrel\Psi\lra\; (W\otimes \eps^n)\gotimes\CO(n,0),
$$
\begin{equation}\label{psiphieq}
\Phi = (0\;,\;-\widehat{(B_1z-x)}I\;,\;P(x,z))\;,\qquad
\Psi = (-zwJ\widehat{(B_1z-x)}\;,\;0\;,\;P(x,z)),
\end{equation}
where $\widehat{(B_1z-x)}$ stands for the cofactor 
 matrix (i.e., the matrix formed by the $(n-1)\times (n-1)$ minors 
in the matrix $B_1z-x,$ taken with appropriate sign)
and $P(x,z)=\det(B_1z-x)$. It is easy to see that
$\Psi\cdot a = b\cdot \Phi = 0$. Thus $\Phi$ and $\Psi$
induce morphisms
$$
(W\otimes \eps^{-n})\gotimes\CO(-n,0) \stackrel\phi\lra E
\stackrel\psi\lra (W\otimes \eps^n)\gotimes\CO(n,0).
$$
Furthermore, it is easy to show that the composite
$\psi\ccirc\phi:\,(W\otimes \eps^{-n})\gotimes\CO(-n,0)\to
(W\otimes \eps^n)\gotimes\CO(n,0)$ equals multiplication
by $P(x,z)^2$.
Thus $(\phi,\psi)$ is a trivialization of $E$ in a neighborhood 
of $\PP^1_z$. Finally, it is easy to see that the trivialization
of $i_w^*E\cong W\gotimes\CO_{\PP^1_w}$ takes form
$$
(W\otimes \eps^{-n})\gotimes\CO(-n) \stackrel{P(x,z)}\lra
W\gotimes\CO \stackrel{P(x,z)}\lra (W\otimes \eps^n)\gotimes\CO(n).
$$
This trivialization is equivalent to the trivial one, thus
the map $i_w^*:\Gradp(W)\to\Graff(W)$ takes $(E,\phi,\psi)$
to the base point $\CW_0$.
Thus we obtain an embedding
$$
\beta:\;\;\bigsqcup\nolimits_V\;
\MPPP(V,W)\; \into\; (\iw^*)^{-1}(\CW_0)\; \subset\; \Gradp(W)\;.
$$

Theorem \ref{framedsheaves} is an immediate consequence of
 Theorem~\ref{qs1}
and  the following result.

\begin{theorem}\label{betaiso}
The map 
$\beta:\;\;\bigsqcup\nolimits_V\;\MPPP(V,W)\; \too\; (\iw^*)^{-1}(\CW_0)$ 
is a bijection.
\end{theorem}
{\sl Proof:}
Let $E$ be a coherent 
sheaf on $\PPP$ with a trivialization $(\phi,\psi)$ in
a neighborhood of $\PP^1_z$. Then~$E$ has a canonical
$W$-framing on $\PP^1_z$ (given by restricting the
trivialization). If, in addition, $\iw^*(E) = \CW_0$
then the sheaf $E$ acquires a canonical $W$-framing
on $\PP^1_w$. Moreover, the framings agree at the point~$P$,
hence we obtain a map
$$
\alpha:\;\; (\iw^*)^{-1}(\CW_0)\; \too\; \bigsqcup\nolimits_V\;\MPPP(V,W).
$$

We now show that both $\alpha\ccirc\beta$ and $\beta\ccirc\alpha$
are identities. To prove that $\alpha\ccirc\beta=\id$ note
that the $W$-framings of the sheaf $E$ on $\PP^1_z$ and $\PP^1_w$
induced by the trivialization~(\ref{psiphieq}) coincide with
the canonical $W$-framings. 

In order to prove $\beta\ccirc\alpha=\id$
one needs to show that any trivialization of $E$ which
gives the canonical $W$-framing of $E$ on $\PP^1_w$
is equivalent to the trivialization~(\ref{psiphieq}).
Indeed, let $E$ be the cohomology of  monad~(\ref{monadppp})
and consider an arbitrary trivialization
$(W\otimes \eps^{-n'})\gotimes\CO(-n',0) \stackrel{\phi'}\lra E
\stackrel{\psi'}\lra (W\otimes \eps^{n'})\gotimes\CO(n',0)$
of $E$ in a neighborhood of $\PP^1_z$.
Applying  diagram (\ref{monadppp}) to compute
$\Hom\bigl((W\otimes \eps^{-n'})\gotimes\CO(-n',0),E\bigr)$ and
$\Hom\bigl(E,(W\otimes \eps^{n'})\gotimes\CO(n',0)\bigr)$
we see that the morphisms $\phi'$ and $\psi'$ can be lifted
to morphisms
$$
(W\otimes \eps^{-n'})\gotimes\CO(-n',0)\; \stackrel{\Phi'}\lra
\,\begin{array}{c}
(V\otimes \eps)\gotimes\CO(0,-1) \\
\bigoplus \\
(V\otimes \eps^{-1})\gotimes\CO(-1,0) \\
\bigoplus \\
W\gotimes\CO
\end{array}\,
\stackrel{\Psi'}\lra\; (W\otimes \eps^n)\gotimes\CO(n',0)
$$
such that
\begin{equation}\label{abp}
b\cdot\Phi'=\Psi'\cdot a=0.
\end{equation}
Moreover, the lift $\Phi'$ of $\phi$ is unique, while the lift
$\Psi'$ of $\psi'$ is  unique up to a summand of the form
$\lambda\cdot b,$
where  $\lambda\in\Hom\bigl(V\gotimes\CO,(W\otimes
\eps^{n'})\gotimes\CO(n'
,0)\bigr)$.
Let
$$
\Phi' = (\Phi'_1,\Phi'_2,\Phi'_3),\qquad
\Psi' = (\Psi'_1,\Psi'_2,\Psi'_3)
$$
be the components of  $\Phi'$ and $\Psi'$ with respect
to the direct sum decomposition 
$$(V\otimes \eps)\gotimes\CO(0,-1)\; \bigoplus\;
(V\otimes \eps^{-1})\gotimes\CO(-1,0)\; \bigoplus \;W\gotimes\CO\;.
$$
Since the trivialization $(\phi',\psi')$ of $E$
restricts to the canonical $W$-framing of $E_{\PP^1_w}$
it follows that
$$
\Phi'_3 = \Psi'_3 = P'(x,z)\;\;,
$$
where  $P'(x,z)$ is a  homogeneous polynomial
of certain degree $n'$,
such that $P'(1,0)=1$. Furthermore, the vanishing:
$$
\Hom\bigl((W\otimes \eps^{-n'})\gotimes\CO(-n',0)\,,\,(V\otimes
\eps^{-1})\gotimes\CO(0,-1)\bigr)
= 0
$$
implies  that $\Phi'_1=0$. On the other hand,
using the freedom in the choice of $\lambda$ one
can make $\Psi'_2=0$. Then, equations~(\ref{abp}) yield
$$
(B_1z-x)\Phi'_2 + I P'(x,z) = 0,\qquad
\Psi'_1(B_1z-x) + P'(x,z)zw J = 0.
$$
Multiplying the first equation by $\widehat{(B_1z-x)}$ on
the left and the second by $\widehat{(B_1z-x)}$ on the right
we obtain
$$
P(x,z)\Phi'_2 = - P'(x,z)\widehat{(B_1z-x)} I,\qquad
\Psi'_1 P(x,z) = - zw J\widehat{(B_1z-x)} P'(x,z),
$$
where $P(x,z)=\det(B_1z-x)$. Thus, we see that
the trivialization $(\phi',\psi')$ is equivalent
to trivialization~(\ref{psiphieq}) via the equivalence
given by the polynomials $P(x,z)$ and $P'(x,z)$.
\hfill\qed

\section{Projective $\B$-modules.}

In this section we prove Theorem \ref{prb}.
Throughout, we will assume the parameter $\tau$ to be generic.
We begin with a description of projective $\Brat$-modules.

\begin{proposition}\label{prbrat}
Assume that $\tau$ is generic. Then

\vi Any projective finitely generated $\Brat$-module has the form:
 $\CM\cong W\gotimes\Brat,$ for a  finite
dimensional $\Gamma$-module $W$. 

\vii Two
$\Brat$-modules $W\gotimes\Brat$ and $W'\gotimes\Brat$
are isomorphic if and only if  $\dim W = \dim W'$.
\end{proposition}
{\sl Proof:}
Let $\e=\frac{1}{|\G|}\sum_{\gamma\in \G}\,\gamma\,
\in\CG$ denote the  averaging idempotent.
 Consider the subalgebra
$\e\Brat \e\subset\Brat$. It is clear that
this algebra is isomorphic to the algebra of
differential operators on $\C/\G$ with rational
coefficients. Set:
$\,\xi=
\e x^m ,$ and $\,
\eta=\e x^{1-m}y.
$
We have an isomorphism
$$
\e\Brat \e \;\cong\;
 \C(\xi)\langle\eta\rangle\big/\langle[\eta,\xi]=|\tau|\rangle
$$
($\xi$ can be considered as a coordinate on $\C/\G$ and $\eta$ as
a vector field on $\C/\G$ generating the algebra of differential
operators). We see that $\e\Brat \e$ is a skew polynomial
ring over the field $\C(\xi)$, hence it is Euclidean. Therefore,
$\e\Brat \e$ is a principal ideal domain, hence any projective
$\e\Brat \e$-module is free.

We claim next that the algebras $\Brat$ and $\e\Brat \e$ are Morita
equivalent. To prove this, observe first that
since $\G$-action on $\C\smallsetminus\{0\}$ is free,
the field $\C(x)$ is a Galois extension of $\C(x)^\G,$
with $\G$ being the Galois group. It follows that
the algebra $\CG(x)=\C(x)\#\G$ is a simple
$\C(x)^\G$-algebra. Hence, $\CG(x)\cdot\e\cdot\CG(x)$,
a two-sided ideal in $\CG(x)$, must be equal to $\CG(x)$.
We see that there exist elements $a_j,b_j\in\CG(x)\,,\,j=1,\ldots,l,$
such that $\sum a_j\cdot\e\cdot b_j=1$.
Therefore, since $a_j,b_j\in\CG(x)\subset \Brat$, we deduce
$\Brat\cdot\e\cdot\Brat=\Brat.$
This implies, by a standard argument, that 
the functor $N\mapsto N\otimes_{_{\e\Brat \e}}\e\Brat$ provides
a Morita
equivalence between  the algebras  $\e\Brat \e$ and  $\Brat$.
Our claim follows.

Using the Morita equivalence we deduce that
any projective $\Brat$-module is isomorphic to
$$
(\e\cd\Brat\cd \e)^{\oplus r} \otimes_{\e\cdot \Brat\cdot  \e} (\e\cd\Brat)\; \cong
\;(\e\cd\Brat)^{\oplus r}\; \cong\;
 ({\mathsf{triv}}^{\oplus r})\gotimes\Brat,
$$
where ${\mathsf{triv}}=\eps^0$ 
is the trivial 1-dimensional $\G$-module.
This proves the first part of the Proposition.

To prove the second part, let
$W\cong \oplus W_i\otimes \eps^i$ be 
 a decomposition of $W$
with respect to the irreducible $\G$-modules $\eps^i$.
Then,
$W\gotimes\Brat$ goes, under  Morita equivalence to
$$
(W\gotimes\Brat)\otimes_{\Brat}\Brat \e \;\cong\;
W\gotimes\Brat \e \;=\;
(\bigoplus\, W_i\otimes \eps^i)\gotimes\Brat \e\; \cong\;
\bigoplus\, W_i\otimes \e_i\Brat \e,
$$
where $\e_i\in\CG$ is the projector onto $\eps^i$. Now it is easy
to see that $\e_i\Brat \e$ is a free rank 1 $\e\Brat \e$-module
(with $\e\cd x^i$ being a generator). Hence
$$
\bigoplus\nolimits_i\, W_i\otimes \e_i\Brat \e\;
\; \cong\;\; (\e\Brat \e)^{\oplus\dim W}.
$$
In particular, it follows that $W\gotimes\Brat$ and $W'\gotimes\Brat$ are
isomorphic $\Brat$-modules if and only if  
$\dim W=\dim W'.$~\hfill\qed
\medskip

Recall that 
we have a natural isomorphism $K(\CG)\iso K(\B)\,,\,
W\mapsto W\gotimes\B$.
Let $[N]\in K(\CG)$ denote the class of
a $\B$-module $N$ under the inverse isomorphism, and write
 $\dim:K(\CG)\to\Z$
for the dimension homomorphism.
Proposition~\ref{prbrat} can be reformulated as follows.

\begin{corollary}\label{kbrat}
We have $K(\Brat)=\Z$. Moreover, the morphism $K(\B)\to K(\Brat)$
induced by the localization functor $N\mapsto N\otimes_\B\Brat$
gets identified with the dimension homomorphism $\dim:K(\CG)\to\Z$.
\end{corollary}

\begin{lemma}\label{nrat}
If $N$ is a projective finitely generated $\B$-module
such that $[N]=W$, then $N\otimes_\B\Brat\cong W\gotimes\Brat$.
\end{lemma}
{\sl Proof:}
It is clear that $N\otimes_\B\Brat$ is a projective finitely 
generated $\Brat$-module, hence Proposition~\ref{prbrat} (i) 
yields: $N\otimes_\B\Brat\cong W'\gotimes\Brat$ for 
some~$W'$. Moreover, Corollary~\ref{kbrat} implies that 
$\dim W'=\dim W$. Finally, Proposition~\ref{prbrat} (ii) 
shows that $N\otimes_\B\Brat\cong W\gotimes\Brat$.
\hfill\qed\medskip

\begin{lemma}\label{prfat}
Let $W$ be a $\Gamma$-module. Any projective finitely generated
$\B$-module $N$ with $\dim[N]=\dim W$ can be embedded into
$W\gotimes\Brat$,
as a fat
$\B$-submodule. Furthermore, the embedding
is unique up to the action of the group $\Gw$.
\end{lemma}
{\sl Proof:}
To prove the existence of embedding we consider the natural
map $N \to N\otimes_{\B}\Brat = N\otimes_{\CG[x]}\CG(x)$.
Since $N$ is projective it follows that $N$ is torsion free
(as $\CG[x]$-module), hence the above map is an embedding.
Further, by Lemma \ref{nrat} it follows that
$N\otimes_\B\Brat\cong W\gotimes\Brat$. Finally, taking
an arbitrary set of generators (over $\B$) of $N\subset W\gotimes\Brat$ 
and denoting by $p_1(x)$ some $\Gamma$-invariant multiple
of all their denominators we see that
$$
N \subset \mbox{$\frac{1}{p_1}$}\cd(W\gotimes \B)
\subset W\gotimes\Brat.
$$
Similarly, considering the dual $\B$-module one can check
that there exists a $\Gamma$-invariant polynomial $p_2(x)$
such that 
$$
p_2\cd(W\gotimes\B) \subset N \subset W\gotimes\Brat.
$$
Finally, taking $p(x)=p_1(x)p_2(x)$ we see that $N$ is 
a fat $\B$-submodule in $W\gotimes\Brat$.

Now assume that we have two embeddings
$\psi_1,\psi_2:\,N \into W\gotimes\Brat$. Tensoring with $\Brat$
we obtain two isomorphisms
$\psi_1,\psi_2:\,N\otimes_\B\Brat \iso W\gotimes\Brat$.
Then $g=\psi_2\ccirc\psi_1^{-1}\in GL_\Brat(W\gotimes\Brat)=\Gw$
and it is clear that $\psi_2 = g\ccirc\psi_1$.
\hfill\qed\bigskip

\noindent
{\sl Proof of Theorem~$\ref{prb}$:}
It follows from Lemma~\ref{prfat} that any projective $\B$-module $N$
such that $\dim[N]=r$ can be embedded into $W\gotimes\Brat$ as
a fat $\B$-submodule. On the other hand, for generic $\tau$
the homological dimension of the algebra $\B$ equals~1 (see~\cite{CBH}),
hence any fat $\B$-submodule $N\subset W\gotimes\Brat$ is projective.
Moreover, it is clear that we have 
$N\otimes_\B\Brat=W\gotimes\Brat$, 
hence by Proposition~\ref{prbrat} and Corollary~\ref{kbrat} 
we deduce: $\dim [N]=\dim W=r$. Finally, by Lemma~\ref{prfat}
two fat $\B$-submodules in $W\gotimes\Brat$ are isomorphic
as $\B$-modules if and only if they are conjugate by
 the action of the group $\Gw$. It follows that
the set of isomorphism classes of projective $\B$-modules
$N$ with $\dim[N]=r$ is in a natural bijection with the coset space
$\Gw\backslash\Gradb(W)$. It remains to apply the isomorphisms
of diagram \eqref{gwequi}.
\hfill\qed

\section{Appendix A: Formalism of Polygraded Algebras}

Let $A = \oplus_{p\ge 0} A_p$ be a graded algebra over a field $\k$.
Let $\gr(A)$ denote the category of graded finitely generated
$A$-modules. For any $n\in\Z$ and any $M\in\gr(A)$ let
$M_{\ge n} = \oplus_{p\ge n}M_p$ be the tail of $M$.
An element $x\in M$ is called torsion if $x\cdot A_{\ge n}=0,$
for some $n$. A module $M$ is called torsion if every element 
of $M$ is torsion. Let $\tor(A)$ denote the full subcategory
of $\gr(A)$ formed by all torsion $A$-modules. Then $\tor(A)$
is a dense subcategory, hence one can consider the quotient
category $\qgr(A) = \gr(A)/\tor(A)$. If $A$ is commutative
and generated over $A_0$ by $A_1$ then by the Serre Theorem
the category $\qgr(A)$ is equivalent to the category $\coh(X)$
of coherent sheaves on  $X=\Proj(A)$,
the projective spectrum
of the algebra $A$.

In the series of papers \cite{AZ}, \cite{YZ}, etc., 
 a formalism has been
developped  that allows to consider the category
$\qgr(A)$ as a category of coherent sheaves in the case when
$A$ is a noncommutative graded algebra. This means
 that the category $\qgr(A)$ shares
many of the general properties of categories of coherent sheaves,
provided  the algebra $A$ satisfies some `reasonable' properties.
In this case one says that $\qgr(A)$ is the category of coherent
sheaves on a noncommutative algebraic variety $X$ and denotes
it by $\coh(X)$.

\nc{\bo}{{\mathbf{0}}}
\nc{\bd}{{\mathbf{d}}}
\nc{\bl}{{\mathbf{l}}}
\nc{\bn}{{\mathbf{n}}}
\nc{\bp}{{\mathbf{p}}}
\nc{\bq}{{\mathbf{q}}}
\nc{\N}{{\mathbb{N}}}
\nc{\Id}{{\mathbf{Id}}}

We extend the formalism of \cite{AZ} to the poly-graded case
as follows. Let
$A = \oplus_{\bp\in\N^r} A_\bp$ be an $\N^r$-graded
algebra (we will denote vector indices by bold letters).
Let $\gr^r(A)$ denote the category of finitely generated $\Z^r$-graded
$A$-modules. For any $\bn\in\Z^r$ and any $M\in\gr^r(A)$
let $M_{\ge\bn} = \oplus_{\bp\ge\bn}M_\bp$ be the tail of $M$,
where $\bp=(\bp_1,\dots,\bp_r)\ge\bn=(\bn_1,\dots,\bn_r)$
if and only if  $\bp_i\ge\bn_i$ for all $1\le i\le r$. An element
$x\in M$ is called torsion if $x\cdot A_{\ge\bn}=0,$
for some $\bn$. A module $M$ is called torsion if every its
element is torsion. Let $\tor^r(A)$ denote the full subcategory
of $\gr^r(A)$ formed by all torsion $A$-modules. Thus, $\tor^r(A)$
is a fat subcategory. We let
$\qgr^r(A) = \gr^r(A)/\tor^r(A)$ be a polygraded counterpart
of the category.

In the polygraded situation we have to make
 the following modifications in the definitions used in  \cite{AZ}.
First, a $\Z^r$-graded $\k$-module $V=\oplus_{\bp\in\Z^r}V_\bp$
should be called left bounded if $V = V_{\ge\bn}$ for some $\bn\in\Z^r$
(such $\bn$ is called a left bound for $V$). Similarly, $V$ 
should be  called
right bounded if $V_{\ge\bn}=0$ for some $\bn\in\Z^r$ (such $\bn$
is called a right bound for $V$). Note, that a finitely generated
module $M$ over a finitely generated algebra $A$ is torsion if and only if 
it is both left and right bounded. Thus $\tor^r(A)$ is the category
of bounded $\Z^r$-graded $A$-modules.

Most essential  changes involve the definition
of {\sf property} $\bchi_i(M)$, cf. (\cite{AZ}, Definition~3.2).
First, introduce the following notation. For each $i=1,\dots,r,$
write  $e_i\in\Z^r$ for  the $i$-th basis vector, and let
 $I\subset\{1,\dots,r\}$
 denote
a nonempty subset of indices. For any $M\in\gr^r(A)$ put
$$
M_\bn^I =
\Big(\bigoplus_{\bp\ge\bn,\ \bp_i=\bn_i\ \text{for $i\in I$}} M_\bp\Big) =
M_{\ge\bn}\Big/\sum\nolimits_{i\in I}\;M_{\ge\bn+e_i}\;\;.
$$

\begin{definition}
We say that  property $\bchi_i(M)$ holds for a $\Z^r$-graded $A$-module $M$
provided\linebreak
 $\Ext^j(A_\bo^{\{k\}},M)$ is bounded for all $j\le i$ and
all $1\le k\le r$. 

We say that  property  $\bchi_i$
holds for the graded algebra $A$
provided   property $\bchi_i(M)$ holds for every finitely
generated $\Z^r$-graded $A$-module $M$.

We say that  property  $\bchi$ holds for $A$ provided
 property $\bchi_i$ holds
for every $i$.
\end{definition}

In \cite{AZ},  a graded  algebra $A$ was said to be
 {\it regular of dimension $d$}
 if the following holds
\smallskip

\noindent$(0)$
$A$ is connected (i.e.\ $A_0=\k$);

\noindent$(1)$
$A$ has finite global dimension $d$;

\noindent$(2)$
$A$ has polynomial growth;

\noindent$(3)$
$A$ is Gorenstein, that is
$\Ext^i_{\mmod(A)}(\k,A) = \begin{case}
\k[l], & \text{if $i=d$}\\
0, & \text{otherwise}
\end{case}$
\medskip

It was demonstrated in \cite{AZ} that, for regular algebras $A$,
the category $\qgr(A)$ has good properties,
in particular, one can
 compute cohomology of the sheaves $\CO(i) = \pi(A(i))$, where
$\pi:\gr(A)\to\qgr(A)$ is the projection functor and $(i)$ stands
for the degree-shift  by $i\in\Z$. Further,  in \cite{BGK}
we explained that one can replace conditions (0) and (3) above
 by the following conditions:

\medskip
\noindent$(0')$
$A_0$ is a finite dimensional semisimple $\k$-algebra;

\noindent$(3')$
$A$ is generalized Gorenstein, that is
$\Ext^i_{\mmod(A)}(\k,A) = \begin{case}
R[l], & \text{if $i=d$}\\
0, & \text{otherwise}
\end{case}$

\medskip
\noindent
where $R$ is a finite dimensional $A_0$-bimodule
isomorphic to $A_0$ as right $A_0$-module.

\medskip

In this paper we will need a further generalization
of the notion of regular algebra to the setup of polygraded
algebras. To this end,  one 
has to replace  condition $(3')$ above
by the following
condition:

\medskip
\noindent$(3'')$
$A$ is strongly Gorenstein with parameters
$\bd = (d_1,\dots,d_r)$, $\bl = (l_1,\dots,l_r)$,
such that $d = \sum_{i=1}^r d_i$. This means that
for any subset $I\subset\{1,\dots,r\}$ we have
$$
\Ext^i_{\mmod(A)}(A_\bo^{I},A) = \begin{case}
(R_I\otimes_{A_0}A_\bo^{I})(\bl_I), & \text{if $i=d_I$}\\
0, & \text{otherwise}
\end{case}
$$
where 
$$
d_I = \sum\nolimits_{k\in I}\;d_k\quad,\quad
\bl_I = \sum\nolimits_{k\in I}\;l_ke_k\quad,\quad
R_I = \bigotimes\nolimits_{_{k\in I}}\;\eps^k\quad
(\text{tensor product over}\enspace A_0),
$$
%where $d_k$ and $l_k$ are integers, and
where $\eps^k$ are $A_0$-bimodules isomorphic to $A_0$ as right $A_0$-modules,
and 
such that $\eps^k\otimes_{A_0}\eps^l$ $\cong \eps^l\otimes_{A_0}\eps^k$ as $A_0$-bimodules.

Now, with all these modifications made, one can verify that most of the
results of \cite{AZ} can be extended to
$\N^r$-graded algebras by the same arguments as in  \cite{AZ}.
In particular, we have an analog of (\cite{AZ}, Theorem 8.1).

\begin{theorem}\label{afteraz}
Let $A$ be a $\N^r$-graded noetherian regular algebra of dimension $d$
over a semisimple algebra $A_\bo$.
Let $\CO(\bp)=\pi(A(\bp))\in\qgr^r(A) = \coh(X)$.
Then

\noindent$(1)$ Property  $\bchi$ holds for $A$.

\noindent$(2)$ $H^0(X,\CO(\bp)) = A_\bp$ and
$H^{>0}(X,\CO(\bp))=0,$ for all $\bp\ge0$.

\noindent$(3)$ The cohomological dimension of the category
$\coh(X) = \qgr^r(A)$ equals $d-r$.\hfill\qed
\end{theorem}

\begin{remark} As opposed to the single-graded case
studied in \cite{AZ}, in the polygraded case
it is  impossible  to
determine the cohomology of the sheaves $\CO(\bp)$ for nonpositive
$\bp$ without some
extra information about the structure of the algebra $A$
(it is necessary to know the $A_\bo^I$-module structure
on $A_\bn^I$ for all $\bn\ge0$ and $I\subset\{1,\dots,r\}$).
\end{remark}

\begin{definition}
We say that an $\N^r$-graded algebra $A$ is strongly generated
by its first component if for any $1\le i\le r$ both maps
below are surjective for any $\bp\ge\bo$
$$
A_{e_i}\otimes A_\bp \too A_{\bp+e_i}
\quad\text{and}\quad
A_\bp\otimes A_{e_i} \too A_{\bp+e_i}\;\;.
$$
\end{definition}

\begin{remark}
It is easy to see that any $\N$-graded algebra which is
generated
by its first component is strongly generated. Thus in the
case $r=1$ we obtain nothing new.
\end{remark}

An element $\bp=(p_1,\dots,p_r)\in\N^r$ is called strictly
positive if $p_i>0$ for all $1\leq r$.

\begin{proposition}\label{bpample}
If $A$ is an $\N^r$-graded noetherian algebra strongly generated
by its first component and satisfying the  $\bchi$-condition, then
for any strictly positive $\bp$ the shift functor $s(M)=M(\bp)$
in the category $\qgr^r(A)$ is ample in the sense of
{\rm\cite{AZ}, (4.2.1)}.
\end{proposition}
{\sl Proof:}
It follows from an $\N^r$-graded analog of the
Theorem 4.5 of loc.cit.\ that the collection of shift functors
$s_i(M)=M(e_i)$, $i=1,\dots,r$ is ample.
Now let $\CE$ be an object of $\qgr^r(A)$. Then it follows
from the ampleness of the collection $(s_i)$ that
there exists a surjection $\oplus_{i=1}^p\CO(-\bl_i)\to\CE$
for some $\bl_i\ge\bo$. Now for each $\bl_i$ we can choose
$k_i\in\N$ such that $k_i\cdot\bp\ge\bl_i$. Then the strong 
regularity of the algebra $A$ implies that the canonical map
$$
A_{k_i\cdot\bp-\bl_i}\otimes_{A_0}\CO(-k_i\cdot\bp) \too \CO(-\bl_i)
$$
is surjective. Further, since $A_{k_i\cdot\bp-\bl_i}$ is
a finitely generated $A_0$-module it follows that we have
a surjection $\oplus_{i=1}^p\CO(-k_i\cdot\bp)^{\oplus m_i}\to\CE$
and  part (a) of the ampleness property for the functor $s$ follows.

Part (b) of the ampleness for the functor $s$ follows trivially
from the ampleness of the collection $s_i$.
\hfill\qed\medskip

\begin{remark} For any strictly positive $\bp$ we put
$\,
\Delta_\bp(A) := \bigoplus_{k=0}^\infty\; A_{k\cdot\bp}
.\,$
Thus, $\Delta_\bp(A)$ is a {\it single-graded}
 subalgebra of $A$.
The following is immediate from Proposition \ref{bpample} and
(\cite{AZ}, Theorem~4.5).
\end{remark}

\begin{corollary}\label{qgrrqgr}
If $A$ is an $\N^r$-graded noetherian algebra strongly generated
by its first component and satisfying the condition $\bchi$, then
for any strictly positive $\bp$ the algebra $\Delta_\bp(A)$ is
noetherian, satisfies the condition $\bchi$ and we have
an equivalence of categories
$$
\qgr^r(A) \cong \qgr(\Delta_\bp(A)).\hskip 30mm\hfill\square
$$
\end{corollary}

\begin{remark}
We would like to emphasize that, inspite of Corollary~\ref{qgrrqgr},
 the   above developed formalism of  quotient
categories for polygraded algebras  {\it does not} reduce to that
for single-graded   algebras. The point is that though  algebras
$A$ and $\Delta_\bp(A)$ give rise to equivalent
  quotient
categories, the  algebra  $\Delta_\bp(A)$ may not be regular or Koszul, for
instance, even when
$A$ is.
\end{remark}

\section{Appendix B: The Geometry of $\PPP$}
The goal of this Appendix is to study the homological properties of the
algebra $\Qu$, see \eqref{qu},
 and to establish Serre Duality and
Beilinson Spectral Sequence for $\PPP$.

\begin{proposition}\label{Qreg}
The bigraded algebra $\Qu$ is noetherian and is strongly generated 
by its first component. Furthermore, $\Qu$ is regular of dimension $4$.
\end{proposition}

To prove this Proposition we introduce some notation.
Given  a semisimple algebra $A_\bo$ and an $A_\bo$-bimodule
$M$, we write $T_{A_\bo}(M)$ for  the tensor algebra
of $M$ over $A_\bo$.

\begin{definition}
Let $A$ be a $\N^r$-graded algebra generated by
$\oplus_{i=1}^r A_{e_i}$ over a semisimple algebra $A_\bo$.
We say that $A$ is quadratic if
$\dis
A = T_{A_\bo}(\oplus_{i=1}^r A_{e_i})\big/\big\langle R\big\rangle,
\,$
where $\big\langle R\big\rangle$ denotes the two-sided
ideal generated by a graded vector subspace
$R = \oplus_{1\le i,j\le r} R_{e_i+e_j}\,$ 
(called `quadratic relations').
\end{definition}

Assume that $A$ is a quadratic $\N^r$-graded algebra.
Let $A^!$ denote its quadratic dual algebra (with respect
to the total grading). Then $A^!$ is also a quadratic
$\N^r$-graded algebra. Recall that the algebra $A$ is
called Koszul if the following  Koszul complex $\K^\bullet(A)$
is exact
$$
\dots \to\,
\oplus_{1\le i,j\le r} (A^!_{e_i+e_j})^*\otimes_{A_0} A(-e_i-e_j) \too
\oplus_{1\le i\le r} (A^!_{e_i})^*\otimes_{A_0} A(-e_i) \too
A \too A_\bo \to 0\;.
$$

\begin{definition}
We call the algebra $A$  strongly Koszul if for
any subset $I\subset\{1,\dots,r\}$ the  following partial
Koszul complex $\K^\bullet_I(A)$ is exact
$$
\dots \to\,
\oplus_{i,j\in I} (A^!_{e_i+e_j})^*\otimes_{A_0} A(-e_i-e_j) \too
\oplus_{i\in I} (A^!_{e_i})^*\otimes_{A_0} A(-e_i) \too
A \too A_\bo^{I} \to 0\;.
$$
\end{definition}

It is clear from the definition of the quadratic dual algebra
that $(A^!)_\bo^{I}$ is dual to $A_\bo^{\bar{I}}$, where
$\bar{I} = \{1,\dots,r\}\setminus I$. Thus if $A$ is strongly
Koszul then for any $I\subset\{1,\dots,r\}$ the algebra
$A_\bo^{I}$ is Koszul as well.
Fix  $\bd=(d_1,\dots,d_r)$, and for any subset $I$,
write  $\bd_I = \sum_{i\in I}d_ie_i$.

\begin{definition}
We say that $A^!$ is strongly Frobenius
of index $\bd$ if the following holds

\vi $A^!_\bp=0$
unless $\bo\le\bp\le\bd$;

\vii The component
$A^!_{\bd_I}$ of $A^!$ is isomorphic to $A_\bo$
as right $A_0$-module, for any
subset $I\subset\{1,\dots,r\}$;

\viii The multiplication map:
$A^!_\bp\otimes_{A_0} A^!_{\bd_I-\bp}\to A^!_{\bd_I}$
gives a nondegenerate pairing, for any $0\le\bp\le\bd_I$.
\end{definition}

\begin{proposition}\label{Qkoszul}
The algebra $\Qu$ is strongly Koszul and
$\Qu^!$ is strongly Frobenius of index $(2,2)$. Moreover, we have
$$
\Qu^!_{i,j} = \begin{case}
\CG, & (i,j)=(0,0)\\
\C\langle\xi,\zeta\rangle\otimes\CG, & (i,j)=(1,0)\\
\C\langle\xi\wedge\zeta\rangle\otimes\CG, & (i,j)=(2,0)\\
\C\langle\eta,\omega\rangle\otimes\CG, & (i,j)=(0,1)\\
\C\langle\xi\wedge\eta,\xi\wedge\omega,
         \zeta\wedge\eta,\zeta\wedge\omega\rangle\otimes\CG, & (i,j)=(1,1)\\
\C\langle\xi\wedge\zeta\wedge\eta,
         \xi\wedge\zeta\wedge\omega\rangle\otimes\CG, & (i,j)=(2,1)\\
\C\langle\eta\wedge\omega\rangle\otimes\CG, & (i,j)=(0,2)\\
\C\langle\xi\wedge\eta\wedge\omega,
         \zeta\wedge\eta\wedge\omega\rangle\otimes\CG, & (i,j)=(1,2)\\
\C\langle\xi\wedge\zeta\wedge\eta\wedge\omega\rangle\otimes\CG, &
(i,j)=(2,2)\\
0, & \text{otherwise}
\end{case}
$$
where $\xi$ is a generator of  $\Gamma$-bimodule $\eps^{-1}$,
$\eta$  is a generator of  $\Gamma$-bimodule $\eps$, and $\zeta,\omega$
are each generators of the trivial $\Gamma$-bimodule.
\end{proposition}
{\sl Proof:}
In the proof we consider the algebra $\Qu$ as an algebra,
depending on a parameter $\tau$. We will indicate the value 
of $\tau$ by a superscript. For example, $\Qu^0$ stands for
the algebra $\Qu$ with $\tau=0$.

First it is easy to show that for any $\tau$ the components
of the dual algebra are given by the above formulas. Further
note that for $\tau = 0$ we have an isomorphism
$\Qu^0 \cong \C[x,z,y,w]\#\Gamma$. In this case it is quite
easy to show that $\Qu^0$ is strongly Koszul. Finally we
note that we may view the family of partial Koszul complexes
$\K^\bullet_I(\Qu^\tau)$ of the algebras $\Qu^\tau$ as a family
of varying (with $\tau$) differentials on the partial Koszul
complex $\K^\bullet_I(\Qu^0)$. Since the complex is exact for
$\tau = 0$ the same is true for all values of $\tau$ close enough
to zero. However, the algebras $\Qu^\tau$ and $\Qu^{\alpha\cdot\tau}$
are isomorphic for any $\alpha \in \C^*$. Thus, $\Qu^\tau$ is
strongly Koszul for any $\tau$.

Similarly, to show that $(\Qu^\tau)^!$ is strongly Frobenius
for any $\tau$ we note that it is true for $\tau=0$.
Further, we consider the family of pairings
$(\Qu^\tau)^!_\bp\gotimes(\Qu^\tau)^!_{\bd_I-\bp}\to(\Qu^\tau)^!_{\bd_I}$
as a family of varying (with $\tau$) pairings
$(\Qu^0)^!_\bp\gotimes(\Qu^0)^!_{\bd_I-\bp}\to(\Qu^0)^!_{\bd_I}$.
Since the pairings are nondegenerate for $\tau=0$ the same
is true for all values of $\tau$ close enough to zero.
However, the algebras $(\Qu^\tau)^!$ and $(\Qu^{\alpha\cdot\tau})^!$
are isomorphic for any $\alpha \in \C^*$. Thus, $(\Qu^\tau)^!$ is
strongly Frobenius for any $\tau$.
\hfill\qed\medskip

\begin{proposition}\label{koszulreg}
If an $\N^r$-graded algebra $A$ is strongly Koszul and
the dual algebra $A^!$ is strongly Frobenius of index $(d_1,\dots,d_r)$
then $A$ is strongly Gorenstein with parameters $\bd=(d_1,\dots,d_r)$
and $\bl=(d_1,\dots,d_r)$.
\end{proposition}
{\sl Proof:}
If $A$ is strongly Koszul then the partial Koszul complex
$\K^\bullet_I(A)$ can be considered as a projective resolution
of $A_\bo^I$. It follows that $\Ext^\bullet_{\mmod(A)}(A_\bo^I,A)$
coincides with the cohomology of complex
$$
0 \to A \,\to\,
\oplus_{i\in I} A^!_{e_i}\otimes_{A_0} A(-e_i) \,\to \dots \to\,
%\oplus_{i,j\in I} A^!_{e_i+e_j}\otimes_{A_0} A(-e_i-e_j) \,\to\, \dots
\oplus_{i\in I} A^!_{\bd_I-e_i}\otimes_{A_0} A(\bd_I-e_i) \,\to\,
A^!_{\bd_I}\otimes_{A_0} A(\bd_I) \to 0.
$$
On the other hand, the strong Frobenius property of the algebra $A^!$
shows that $A^!_\bp\cong A^!_{\bd_I}\otimes_{A_0}(A^!_{\bd_I-\bp})^*$
as $A_0$-bimodule. Hence the above complex is isomorphic to the
complex $A^!_{\bd_I}\otimes_{A_0}\K^\bullet_I(A)(\bd_I)$ truncated 
at the rightmost term.
Therefore, it has a single nonzero cohomology group in  degree $d_I$,
which is isomorphic to $A^!_{\bd_I}\otimes_{A_0}A_\bo^I(\bd_I)$.
It follows that $A$ satisfies the strong Gorenstein property
with parameters $(\bd,\bd)$ and with $R_I=A^!_{\bd_I}$.
\hfill\qed\medskip

{\sl Proof of Proposition~$\ref{Qreg}$:}
It is clear that $\Qu$ is strongly generated by its first component.
So, it remains to prove regularity and the noetherian property.

First, note that $\Qu^\tau_{(0,0)} = \CG$ is a semisimple algebra.
Thus $(0')$ holds.
Second, we have to show that $\Qu^\tau$ is noetherian.
This follows from the fact that $\Qu^\tau$ can
be represented as a consecutive Ore extension of the
base field $\C$.
Further, it is easy to show that
$\dim_\C \Qu^\tau_{i,j} = (i+1)(j+1)|\Gamma|$.
In particular, $\Qu^\tau$ has polynomial growth.
Thus $(2)$ holds.

The strong Gorenstein property $(3'')$ for the algebra
$\Qu^\tau$ follows immediately from Proposition
\ref{Qkoszul} and Proposition
\ref{koszulreg}. The Gorenstein parameters are given by:
$\bd=(2,2)$ and $\bl=(2,2)$.

Finally, it follows from \cite{Hu} that the global dimension
of $\Qu^\tau$ equals the length of the minimal free resolution
of $\Qu^\tau_{(0,0)}$. But the Koszul complex $\K^\bullet(\Qu^\tau)$
provides such resolution of length $4$, hence the global dimension
of $\Qu^\tau$ is lbounded by $4$ from above.
 On the other hand, since $\Qu^\tau$
is Gorenstein with parameters $\bd=(2,2)$, $\bl=(2,2)$ it follows
that $\Ext^4(\Qu^\tau_{0,0},\Qu^\tau)\ne0$, hence the global
dimension equals $4$.
\hfill\qed\medskip

Thus, the cohomological dimension of the category
$\coh(\PPP)=\qgr^2(\Qu^\tau)$ equals $2$, and
it is clear that we have
$$
H^p(\PPP,\CO(i,j)) = \begin{case}\Qu_{i,j}, & p=0\\0, & p>0\end{case}
\qquad\text{for all $i,j\ge0$.}
$$
More generally, we prove
\begin{lemma}\label{hpoij}
$$
H^p(\PPP,\CO(i,j)) = \begin{case}
\Qu_{i,j}, & \text{if $p=0$ and $i,j\ge0$}\\
\eps^{-1}\otimes \Qu_{-2-i,0}^*\gotimes \Qu_{0,j}, &
\text{if $p=1$ and $i\le -2$, $j\ge0$}\\
\eps\otimes \Qu_{0,-2-j}^*\gotimes \Qu_{i,0}, &
\text{if $p=1$ and $i\ge0$, $j\le -2$}\\
\Qu_{-2-i,-2-j}^*, & \text{if $p=2$ and $i,j\le -2$}\\
0, & \text{otherwise}
\end{case}
$$
\end{lemma}
{\sl Sketch of Proof.\;} 
In order to compute the global cohomology of $\CO(i,j)$ for
not necessarily positive values of $(i,j)$ we  use  partial
Koszul complexes. In more detail, the projections to
 the category $\coh(\PPP)$
of the partial
Koszul complexes yield exact sequences
$$
0 \to \eps\otimes\CO(-2,0) \too
\Qu_{1,0}\gotimes\CO(-1,0) \too \CO \to 0,
$$
$$
0 \to \eps^{-1}\otimes\CO(0,-2) \too
\Qu_{0,1}\gotimes\CO(0,-1) \too \CO \to 0
$$
(we used here the fact that $\Qu_\bo^I\in\tor(\Qu)$ for
any nonempty $I\subset\{1,2\}$, and that
$\Qu^!_{0,1}\cong \Qu_{0,1}^*$,
$\Qu^!_{1,0}\cong \Qu_{1,0}^*$,
$\Qu^!_{0,2}\cong \eps^{-1}$,
$\Qu^!_{2,0}\cong \eps$).

To complete the proof of the Lemma we apply descending induction 
in $(i,j)$ using the above sequences twisted by $(i+2,j)$ and 
$(i,j+2)$ respectively, and the fact that the multiplication 
map $\Qu_{i,0}\gotimes \Qu_{0,j}\to \Qu_{i,j}$ is an 
isomorphism of $\Gamma$-bimodules.~\qed
\bigskip

\noindent
{\bf Serre Duality for  $\PPP$.}\;\;
A natural approach to
Serre Duality theorems for
 noncommutative schemes corresponding to regular
noncommutative algebras would be  via the concept of 
balanced dualizing complex (see \cite{Y,YZ}).
Generalizing  this concept to the case of $\N^r$-graded
algebras does not seem  to be straightforward however.
The reason is that,
while the notion of  dualizing complex easily extends to the
polygraded case,  it is
rather difficult to find the relevant meaning of 
 ``balanced'' in this case.
 The problem is
similar to that of computing the cohomology of sheaves $\CO(\bp)$
for nonpositive values of $\bp$, see Remark following Theorem~\ref{afteraz}.

In the special case of $\PPP$ these problems can be circumvented
as follows. We consider the algebra $A=\Delta_{(1,1)}(\Qu)$. It 
follows from Proposition~\ref{Qreg} and Corollary~\ref{qgrrqgr} 
that this algebra is noetherian and satisfies condition $\bchi$. 
Moreover, by  Corollary~\ref{qgrrqgr}, Theorem~\ref{afteraz} and 
Proposition~\ref{Qreg} the cohomological degree of the category 
$\qgr^2(A)=\coh(\PPP)$ equals~2. Hence we can use (\cite{YZ}, 
Theorem~2.3) which implies that the category $\coh(\PPP)$ enjoys 
the Serre duality with dualizing sheaf defined by
$$
\omega^0 = \pi\Big(\oplus_{k=0}^\infty H^2(\PPP,\CO(-k,-k))^*\Big).
$$
But  Lemma \ref{hpoij} yields
$$
\pi\Big(\bigoplus_{k=0}^\infty H^2(\PPP,\CO(-k,-k))^*\Big) \cong
\pi\Big(\bigoplus_{k=0}^\infty \Qu_{k-2,k-2}\Big) \cong
\CO(-2,-2).
$$
Thus the dualizing sheaf on $\PPP$ is isomorphic to $\CO(-2,-2)$.
\bigskip

\noindent
{\bf Beilinson Spectral Sequence.}\;\;
In the noncommutative setting, 
an analogue of Beilinson Spectral Sequence has been
introduced in~\cite{KKO} for  a certain
class of  graded Koszul  algebras, using a double Koszul
bicomplex. Below, we explain how to adapt the approach of \cite{KKO}
to the case of
 $\N^r$-graded  Koszul  algebras. We will freely use the
notation and definitions of  \cite{KKO}, in particular,
the notion of  Yang-Baxter operator.

An exact  functor from  a tensor category $T$
 to the tensor category of vector spaces will be called
a `noncommutative' fiber functor  if this functor is
compatible with the tensor product structures
 and associativity constraint,
but is not necessarily
compatible with the commutativity constraint.
Given a Yang-Baxter operator on a finite-dimensional
vector space, one can construct
as has been explained in \cite{Lyu} (see also \cite{KKO}, section 8),
 a tensor category $T$ equipped with a
`noncommutative' fiber functor. Then the
category of (either graded, or $\N^r$-graded, or \dots) comutative algebras
in the category $T$ gives a class of (graded, $\N^r$-graded, \dots)
noncommutative algebras in the category of vector spaces.
The  class of  noncommutative algebras thus obtained
shares a lot of properties of the category
of commutative algebras. For example, for any two algebras in the
class their tensor product admits a canonical algebra structure.

\begin{remark}
Instead of Yang-Baxter operator in a vector space one may
start with a $A_0$-invariant Yang-Baxter operator in a finitely
generated $A_0$-bimodule, for any semisimple finite dimensional
algebra $A_0$. Then we obtain an $A_0$-linear tensor category
$T$ with a functor to the category of $A_0$-bimodules.
The category of commutative $A_0$-algebras in $T$ gives
a class of noncommutative $A_0$-algebras.
\end{remark}
\medskip

\begin{example}
Consider a free right $\Gamma$-module $V$ of rank $4$ with
generators $x,y,z,w$ and endow it with a $\Gamma$-bimodule
structure as in $(\ref{gxyzw})$. Then the $\Gamma$-linear
operator $V\gotimes V\to V\gotimes V$ defined on the
generators as
$$
\begin{array}{l}
\displaystyle
x\otimes y\mapsto y\otimes x - \frac\tau2z\otimes w - \frac\tau2w\otimes z,
\medskip\\\displaystyle
y\otimes x\mapsto x\otimes y + \frac\tau2z\otimes w + \frac\tau2w\otimes z,
\medskip\\
u\otimes v\mapsto v\otimes u\quad\text{otherwise}
\end{array}
$$
is a Yang-Baxter operator. It is easy to see that the algebra $\Qu$
comes from a bigraded commutative algebra in the tensor category
corresponding to this Yang-Baxter operator.
\end{example}

Let $A$ be an $\N^r$-graded algebra obtained in such a way.
Then $A\otimes_{A_0}A$ is an $\N^r\oplus\N^r$-graded algebra
and the maps $p_1^*(a)=a\otimes1$, $p_2^*(a)=1\otimes a$ are
homomorphisms of algebras $A\to A\otimes_{A_0}A$.
Let $X$ denote the noncommutative variety, corresponding
to the algebra $A$ and let $X\times X$ denote the noncommutative
variety, corresponding to the algebra $A\otimes_{A_0}A$.
Thus $\qgr^r(A)=\coh(X)$, $\qgr^{2r}(A\otimes_{A_0}A)=\coh(X\times X)$.

Now, if $M$ is a right $\N^r$-graded $A$-module we define
$p_1^*M = M\otimes_A(A\otimes_{A_0}A)$. Then $p_1^*M$ is
a right $(\N^r\oplus\N^r)$-graded $A\otimes_{A_0}A$-module.
It is clear that $p_1^*M(\bp,\bq) = M_\bp\otimes_{A_0}A_\bq$,
hence for any $M\in\tor^r(A)$ we have
$p_1^*M\in\tor^{2r}(A\otimes_{A_0}A)$. Thus $p_1^*$ can be
considered as a functor $\qgr^r(A)\to\qgr^{2r}(A\otimes_{A_0}A)$,
that is, a functor: $\coh(X)\to\coh(X\times X)$.

Similarly, if $M = \oplus_{\bp,\bq}M_{\bp,\bq}$ is a right
$(\N^r\oplus\N^r)$-graded $(A\otimes_{A_0}A)$-module then
we define $((p_2)_*M)_\bq = \Gamma(X,\pi(\oplus_{\bp}M_{\bp,\bq}))$,
where
$$
\Gamma(X,\pi(\bullet)) = \Hom_{\qgr^r(A)}(\pi(A),\pi(\bullet)),
$$
and the $A$-module structure of $\oplus_{\bp}M_{\bp,\bq}$
is obtained from the homomorphism $p_1^*$. It is clear that
$(p_2)_*M = \oplus_\bq((p_2)_*M)_\bq$ is an $\N^r$-graded $A$-module.
Furthermore, if $M\in\tor^{2r}(A\otimes_{A_0}A)$ then it is
clear that $(p_2)_*M\in\tor^r(A)$. Thus $(p_2)_*$ can be
considered as a functor $\qgr^{2r}(A\otimes_{A_0}A)\to\qgr^r(A)$,
that is, a functor: $\coh(X\times X)\to\coh(X)$.

Now, if $N$ is an $(\N^r\oplus\N^r)$-graded $(A\otimes_{A_0}A)$-bimodule
then $M\mapsto (p_2)_*(p_1^*M\otimes_{A\otimes_{A_0}A}N)$ gives a functor
$\Phi_N:\qgr^r(A)\to\qgr^r(A)$, that is, a functor:
 $\coh(X)\to\coh(X)$.

\begin{lemma}\label{PhiN}
\vi Let $(\Delta_A)_{\bp,\bq} = A_{\bp+\bq}$, and
$\Delta_A=\oplus_{\bp,\bq\ge0}(\Delta_A)_{\bp,\bq}$.
Then there is a natural isomorphism
of functors $\Phi_{\Delta_A}\cong\Id$.

\vii
If $N_1$, $N_2$ are $A$-bimodules and $N_1\otimes_{A_0}N_2$
has a canonical structure of an $(A\otimes_{A_0}A)$-bimodule
then
$$
\Phi_{N_1\otimes_{A_0}N_2}(M) =
\Gamma(X,\pi(M\otimes_AN_1))\otimes_{A_0}N_2.
$$
\end{lemma}
{\sl Proof:}
\vi Let $M'=\Phi_{\Delta_A}(M)$.
Note that $\Delta_A$, considered as an $A$-module,
is isomorphic to $\oplus_{\bq\in\N^r}\, A(\bq)_{\ge\bo}$.
Hence $M'_\bq=\Gamma(X,\pi(M\otimes_A A(\bq)_{\ge\bo}))$.
On the other hand, it is clear that
$$
\pi(M\otimes_A A(\bq)_{\ge\bo})\cong
\pi(M\otimes_A A(\bq)) \cong
\pi(M(\bq)) \cong \pi(M)(\bq).
$$
This means that $M' = \oplus_{\bq}\Gamma(X,\pi(M)(\bq))$,
hence $\pi(M')\cong \pi(M)$. Furthermore, it is clear that the
isomorphism that we have constructed gives an isomorphism
of functors $\Phi_{\Delta_A}\to\Id$.

\vii Let $M'=\Phi_{N_1\otimes_{A_0}N_2}(M)$.
Then
$$
M'_\bq=
\Gamma(X,\pi(M\otimes_A N_1\otimes_{A_0}(N_2)_\bq)) =
\Gamma(X,\pi(M\otimes_A N_1))\otimes_{A_0}(N_2)_\bq,
$$
hence $M' = \Gamma(X,\pi(M\otimes_A N_1))\otimes_{A_0}N_2$.
\hfill\qed\medskip

\begin{remark}
It is clear that $\Delta_A$ can be endowed with an
algebra structure. Furthermore, it is easy to show
that $\qgr^{2r}(\Delta_A)\cong\qgr^r(A)$. Finally,
the multiplication in $A$ gives an epimorphism
$A\otimes_{A_0}A\to\Delta_A$. This way, one may
view $\Delta_A$ as a diagonal embedding 
$\Delta_X:X\into X\times X$.
\end{remark}\medskip

Once we have a diagonal  $X \into X\times X$,
we could  apply standard techniques, provided
one finds a resolution of diagonal. If $A$ is Koszul one may
obtain a resolution of diagonal as follows.
Consider the double Koszul bicomplex of $A$:
\nc{\CK}{{\mathcal{K}}}
\nc{\CQ}{{\mathcal{Q}}}
$$
\xymatrix{
{\dots} \ar[r]^-{d_R} &
{\oplus_{i,j} A\otimes (A^!_{e_i+e_j})^*\otimes A(-e_i-e_j) }
\ar[r]^-{d_R} \ar[d]_{d_L} &
{\oplus_i A\otimes (A^!_{e_i})^*\otimes A(-e_i) }
\ar[r]^-{d_R} \ar[d]_{d_L} &
{A\otimes A} \\
\dots \ar[r]^-{d_R} &
{\oplus_{i,j} A(e_i)\otimes (A^!_{e_j})^*\otimes A(-e_i-e_j)}
\ar[r]^-{d_R} \ar[d]_{d_L} &
{\oplus_i A(e_i)\otimes A(-e_i)} \\
\dots \ar[r]^-{d_R} &
{\oplus_{i,j} A(e_i+e_j)\otimes A(-e_i-e_j)}
}
$$
where both $d_R$ and $d_L$ are induced by the differential in the
Koszul complex of $A$. Write
$$
\CK^\bp(A) = \Ker\Big(A(-\bp)\otimes(A^!_\bp)^* \too
\bigoplus\limits_{\{i\ |\ e_i\le\bp\}}
\,
A(e_i-\bp)\otimes(A^!_{\bp-e_i})^*\Big),
$$
for the cohomology of the truncated Koszul complex.
Using the Koszul property of the algebra $A$
and mimicing the proof (\cite{KKO},
Proposition~4.7) we deduce

\begin{proposition}\label{grres}
The following complex is exact
$$
\dots \to
\oplus_{i,j}\,\CK^{e_i+e_j}(e_i+e_j)\otimes A(-e_i-e_j) \to
\oplus_i\,\CK^{e_i}(e_i)\otimes A(-e_i) \to
A\otimes A \to \Delta_A \to 0.
$$
where the map $A\otimes A\to\Delta_A$ is given by the
multiplication in $A$.\hfill\qed
\end{proposition}

Let $\CQ^\bp = \pi(\CK^\bp(A))^*$. Combining \ref{grres}
with \ref{PhiN} we obtain the Beilinson spectral sequence.

\begin{corollary}
Assume that $A$ is Koszul and $A^!$ is Frobenius. Then
for any $F\in X$ there exists a spectral sequence with
the first term
$$
E_1^{p,q} = \bigoplus\limits_{\{\bp\ |\ |\bp|=p\}}
\Ext^q(\CQ^\bp,F)\otimes_{A_0}\CO(-\bp)
\;\;\Longrightarrow\;\;
E_\infty^i = \begin{case}F, & i=0\\0, & \text{otherwise\;.}\end{case}
\hskip 20mm\square$$
\end{corollary}

In the special case of the algebra $\Qu$,
the only non-vanishing components of $\CQ^\bp$ are
$$
\begin{array}{llllll}
\CQ^\bo       & = & \CO, &
\CQ^{e_1}     & = & \eps^{-1}\otimes\CO(1,0),\\
\CQ^{e_2}     & = & \eps^{ 1}\otimes\CO(0,1), \quad &
\CQ^{e_1+e_2} & = & \CO(1,1)
\end{array}
$$
Thus, Beilinson spectral sequence
 takes the form of (\ref{BSSPPP}).

\footnotesize{

}
%\vskip 1cm

\footnotesize{
{\bf V.B.}: Department of Mathematics, California Institute of Technology,
Pasadena CA 91125,
USA;\hfill\break
\hphantom{x}\quad\, {\tt baranovs@caltech.edu}}

\footnotesize{
{\bf V.G.}: Department of Mathematics, University of Chicago,
Chicago IL
60637, USA;\hfill\break
\hphantom{x}\quad\, {\tt ginzburg@math.uchicago.edu}}

\footnotesize{
{\bf A.K.}: Institute for Information Transmission Problems, 
Russian Academy of Sciences,\hfill\break
\hphantom{x}\qquad\qquad\; 19 Bolshoi Karetnyi, Moscow 101447, 
Russia;\quad\, {\tt sasha@kuznetsov.mccme.ru}}

\end{document}